\documentclass[11pt, twoside]{article}
\usepackage{amsfonts,amsmath,amssymb,amsthm}
\usepackage[mathscr]{euscript}
\usepackage{indentfirst}
\usepackage{enumitem}
\usepackage{xcolor}
\usepackage{graphicx}
\usepackage{caption,mathrsfs}
\usepackage{tikz}
\usepackage[hidelinks]{hyperref}

\DeclareMathOperator{\diag}{diag}

\setlist[enumerate]{itemsep=2pt, leftmargin=9mm}
\usepackage{bm}
\usepackage{needspace}
\newtheorem{theo}{Theorem}[section]
\newtheorem{introtheo}{Theorem}[section]
\newtheorem{introcoro}{Corollary}[section]
\newtheorem{prop}{Proposition}[section]
\newtheorem{lem}{Lemma}[section]\newtheorem{lemma}[lem]{Lemma}
\newtheorem{remark}{Remark}[section]
\newtheorem{coro}{Corollary}[section]

\newtheorem{example}{Example}[section]
\numberwithin{equation}{section}

\allowdisplaybreaks[4]

\newcommand{\bea}{\begin{align*}}
\newcommand{\eea}{\end{align*}}
 \newcommand\bee{\begin{equation}}
 \newcommand\eee{\end{equation}}
 \newcommand\bal{\begin{align*}}
 \newcommand\eal{\end{align*}}

\newcommand\ep{\varepsilon}
\newcommand\dd{\displaystyle}

\newcommand\dx{\mathrm{d}x}
\newcommand\dy{\mathrm{d}y}
\newcommand\dt{\mathrm{d}t}
\newcommand\dz{\mathrm{d}z}

\newcommand\qq{\eqref}
\newcommand\kk{\left}
\newcommand\rr{\right}
\newcommand\www{\vspace{-2mm}}
\newcommand\zzz{\vspace{-1mm}}

\newcommand\ud{\underline}
\newcommand\oo{\Omega}\newcommand\boo{\overline\Omega}

\makeatletter
\let\originaltagform@\tagform@
\newcommand{\redtag}{\def\tagform@##1{\originaltagform@{\textcolor{red}{##1}}}}
\newcommand{\redhead}{\thm@headfont{\bfseries\color{red}}}
\makeatother

\usepackage{accents}
\usepackage{amsmath}
\begin{document}
\thispagestyle{empty}
\setlength{\abovedisplayskip}{6.5pt}\setlength{\belowdisplayskip}{6.5pt}
\setlength{\abovedisplayshortskip}{6pt}\setlength{\belowdisplayshortskip}{6pt}

\begin{center}
{\Large\bf Spectral bounds of cooperative nonlocal dispersal systems: monotonicity, domain limits, and applications}\\[3mm]
{ Mingxin Wang$^{a,}$\footnote{{\sl E-mail}: mxwang@sxu.edu.cn. Mingxin Wang was supported by the National Natural Science Foundation of China (NSFC) under Grant 12571221.}
Lei Zhang$^{b,}$\footnote{The corresponding author, {\sl E-mail}: zhanglei890512@gmail.com. Lei Zhang was supported by NSFC under Grants 12471168 and 12471464, the Natural Science Basic Research Program of Shaanxi (2026JC-YXQN-015), and the Fundamental Research Funds for the Central Universities (GK202602001).}\\[2mm]
{\small$^a$ School of Mathematics and Statistics, Shanxi University, Taiyuan 030006, China}\\
{\small$^b$ School of Mathematics and Statistics, Shaanxi Normal University, Xi'an 710119, China}}
\end{center}

\begin{quote}
\noindent{\bf Abstract.} We study the spectral bound of cooperative nonlocal dispersal systems on bounded domains. The spectral bound is continuous under uniform coefficient and suitable nested domain perturbations, and strictly monotone with respect to the coefficients and the domain if either corresponding operator admits a principal eigenvalue. Further results include small-domain limits, large-domain limits for even convolution kernels under suitable coefficient assumptions, and criteria for a unique critical domain size. Applications yield a unique extinction–persistence threshold in the interval length for a heterogeneous epidemic model and a sharp spreading–vanishing threshold for a scalar free-boundary problem. These thresholds remain valid when the spectral bound vanishes on an interval, with the critical habitat length given by its right endpoint.


\end{quote}

\pagestyle{myheadings}
\section{Introduction}
\setlength\arraycolsep{2pt}
\markboth{\hfill Domain dependence of spectral bounds\hfill}
{\hfill M.X. Wang\, \& L. Zhang \hfill}

In spatial ecology and epidemiology, a central question is how habitat geometry and environmental heterogeneity interact with dispersal to govern population persistence or extinction. For reaction-diffusion models driven by the classical Laplacian, this question is answered by the principal eigenvalue of the linearized operator at the trivial equilibrium: because the resolvent is compact, a positive principal eigenfunction always exists, and the principal eigenvalue is strictly monotone with respect to both the reaction coefficients and the domain. Consequently, as a habitat expands, the principal eigenvalue crosses zero at most once, yielding a unique critical domain size that separates extinction from persistence, as well as a sharp spreading--vanishing threshold in free-boundary problems.

When movement is governed by nonlocal dispersal, the differential operator is replaced by an integral operator, and the existence of a principal eigenfunction is no longer guaranteed. Let $n,m\ge1$ and $\mathcal S=\{1,\ldots,m\}$. Given dispersal rates $d_1,\ldots,d_m>0$ with $D=\diag(d_1,\ldots,d_m)$, we consider the cooperative nonlocal dispersal operator $\mathcal J_{A,\oo}$ on a bounded domain $\oo\subset\mathbb R^n$ defined by\vspace{-1mm}
\begin{equation}\label{1.1}
    \begin{cases}
\mathcal J_{A,\oo}[\varphi]=(\mathcal J_{1,A,\oo}[\varphi],\ldots,
\mathcal J_{m,A,\oo}[\varphi])^{\mathsf T},\\[1mm]
\mathcal J_{i,A,\oo}[\varphi](x)
=d_i\displaystyle\int_\oo J_i(x,y)\varphi_i(y)\dy
+A_{i,*}(x)\varphi(x),\;\;\; x\in\boo,
    \end{cases}
\end{equation}
where $A_{i,*}$ denotes the $i$-th row of the matrix $A(x)=(a_{ij}(x))_{m\times m}$. All vectors are column
vectors, and $^{\mathsf T}$ denotes transpose. For $\varphi\in[C(\boo)]^m$, we write $\varphi\gg0$ if $\varphi_i(x)>0$ for every $x\in\boo$ and $i\in\mathcal S$.
We use the same convention on other bounded domains.
Throughout this paper, unless otherwise specified, test functions such as $\varphi$, $\phi$ lie in $[C(\boo)]^m$. Matrix inequalities are understood entrywise.
We impose the following standing hypotheses on the kernels and the coefficient matrix:\zzz
\begin{enumerate}[leftmargin=12mm]
\item[\textbf{(J)}] For each $i\in\mathcal S$, the kernel $J_i:\mathbb R^{2n}\to[0,\infty)$ is continuous and bounded, $J_i(x,x)>0$, and $\int_{\mathbb R^n}J_i(y,x)\dy=1$ for every $x\in\mathbb R^n$.\vspace{-1mm}
\item[\textbf{(A)}] The entries of $A(x)=(a_{ij}(x))_{m\times m}$ belong to $C_b(\mathbb R^n)$. The matrix $A(x)$ is cooperative for every $x$, and $A(x)$ has an irreducibility point in $\mathbb R^n$.\vspace{-1mm}
\end{enumerate}

Unless otherwise specified, all domains considered herein are bounded, with boundaries of zero Lebesgue measure, and contain an irreducibility point of $A(x)$.
For a bounded operator $L$ on $[C(\boo)]^m$, let
$s(L)=\sup\{\operatorname{Re}\lambda:\lambda\in\sigma(L)\}$;
the same notation is used for matrices. For convenience, write
  \[s(A,\overline{\Omega}):=\max_{x\in\boo}s(A(x)).\vspace{-1mm}\]
A real eigenvalue of $\mathcal J_{A,\oo}$ with an eigenfunction $\varphi\gg0$ is called a {\bf principal eigenvalue} and is denoted by $\lambda_p(\mathcal J_{A,\oo})$. Whenever it exists, $\lambda_p(\mathcal J_{A,\oo})=s(\mathcal J_{A,\oo})$.

Although the integral term in \eqref{1.1} is compact on $[C(\boo)]^m$, the multiplication operator associated with $A(x)$ is generally noncompact. Criteria for the existence of principal eigenfunctions have been obtained using generalized Krein--Rutman theory~\cite{EPS72PRSA,N81LNM,BZ07JMAA,C10JDE}, B\"urger's perturbation theory for positive semigroups~\cite{B88MZ,SZ10JDE,SX17DCDSB,BS17PAMS}, and the theory of resolvent-positive operators and their perturbations~\cite{T98DCDS,FLRX24SIMA,LZZ19JDE,KR22MA}. Further results treat partially degenerate systems~\cite{Z24ANS} and cooperative systems on unbounded domains~\cite{WFLS25JDE}. 
When a continuous principal eigenfunction is absent, Coville~\cite{C13AML} studied eigenmeasures.

Variational formulations provide spectral characterizations without assuming the existence of a principal eigenfunction. Results for scalar operators were established in \cite{BCV16JFA,LCW17DCDS}, and variational characterizations for cooperative systems with coupled or uncoupled diffusion were obtained in \cite{SWZ23JDE}. Applications include scalar population equations~\cite{C10JDE,SZ12PAMS}, epidemic models~\cite{YLR19JDE}, stem-cell regeneration~\cite{SLLW23JMPA}, and basic reproduction ratios~\cite{LW24JMAA}. In \cite{WZ26JMPA}, we constructed upper and lower control operators admitting principal eigenfunctions, approximated the generalized principal eigenvalue (i.e., the spectral bound), and established threshold dynamics without requiring a principal eigenvalue for the original operator.

Domain dependence has been studied for several classes of nonlocal operators. Garc\'ia-Meli\'an and Rossi~\cite{GMR09JDE} established strict domain monotonicity, continuity under smooth domain perturbations, and small- and large-domain asymptotics for a scalar nonlocal Dirichlet operator. Berestycki, Coville and Vo~\cite{BCV16JFA} studied domain dependence and dispersal scaling, while Wang and Wang~\cite{WW27AML} investigated dependence on the interval length and coefficients in one dimension. Related results concern spatial heterogeneity and dispersal parameters~\cite{SX15DCDS}, convergence to local diffusion under suitable rescaling~\cite{SX15JDE}, and dispersal limits for cooperative systems~\cite{OSUU24JDE,NV25CVPDE}.

In the absence of a principal eigenfunction, the spectral bound need not be strictly monotone with respect to the coefficients or the domain. It may vanish on an interval of domain sizes, even when both the kernels and the coefficient matrix are symmetric. This phenomenon raises two interconnected questions that motivate the present paper:\www
\begin{enumerate}[label=\textup{(\arabic*)},leftmargin=8mm]
\item  Under what conditions is $s(\mathcal J_{A,\oo})$ strictly monotone with respect to the coefficients and the domain, and what are its small- and large-domain limits?\zzz
\item When the spectral bound vanishes on an interval of habitat lengths,
can the long-term dynamics of the corresponding nonlinear system
still be characterized by a unique critical length?\www
\end{enumerate}

To address Question~(1), we use the control operator approximation from~\cite{WZ26JMPA} and the Collatz--Wielandt formulas. These tools are recalled in Section~\ref{sec2}, together with a dual characterization and a weighted Rayleigh bound. The following theorem provides sufficient conditions for strict coefficient and domain comparisons.\vspace{-1mm}

\begin{introtheo}[Strict comparisons]\label{th1.1} The following statements are valid:\vspace{-2mm}
    \begin{enumerate}[label=\textup{(\roman*)},leftmargin=8mm]
\item \textup{(Coefficient comparison)} Let $B$ satisfy {\bf(A)}, with an irreducibility point in $\oo$, and let $A\ge B$ and $A\not\equiv B$ on $\boo$. If at least one of $s(\mathcal J_{A,\oo})$ and $s(\mathcal J_{B,\oo})$ is the principal eigenvalue, then
  \[ s(\mathcal J_{A,\oo})>s(\mathcal J_{B,\oo}). \]
        
 \item \textup{(Domain comparison)} Let $\oo_1\subset\oo_2$ with $|\oo_2\setminus\oo_1|>0$. If at least one of $s(\mathcal J_{A,\oo_1})$ and $s(\mathcal J_{A,\oo_2})$ is the principal eigenvalue, then
 \[ s(\mathcal J_{A,\oo_1})<s(\mathcal J_{A,\oo_2}).\vspace{-1mm} \]
    \end{enumerate}
\end{introtheo}

If neither operator has a principal eigenvalue, Theorem~\ref{th2.1}\textup{(iii)} implies that each spectral bound equals the corresponding maximum of $s(A(x))$ or $s(B(x))$, and the comparison problem reduces to comparing these maxima.

Lemma~\ref{lem3.2}\textup{(ii)} and Proposition~\ref{prop3.1}
below establish continuity with respect to the coefficients and on suitable nested domains. The corresponding small- and large-domain limits are as follows.\vspace{-1mm}

\begin{introtheo}[Small- and large-domain limits]\label{th1.2} Assume that $x_0$ is an irreducibility point of $A(x)$. Then the following statements are valid:\www
    \begin{enumerate}[label=\textup{(\roman*)},leftmargin=8mm]
\item  $\lim_{\sigma\searrow0}s\big(\mathcal J_{A,B_\sigma(x_0)}\big)=s(A(x_0))$.\zzz        
\item Suppose that $J_i(x,y)=K_i(y-x)$ with $K_i(z)=K_i(-z)$ for every component. Let $A^\infty$ be a constant cooperative irreducible matrix satisfying
\begin{equation}\label{1.2}
  A(x)\le A^\infty,\;\;\;\forall\; x\in\mathbb R^n.
 \end{equation}
  If there exist cubes $H_k=z_k+[0,l_k]^n$ such that
\begin{equation}\label{1.3}
    l_k\to\infty,\;\;\; \max_{i,j\in\mathcal S}\,\sup_{x\in H_k} |a_{ij}(x)-a^\infty_{ij}|\to0,
  \end{equation}
 then
\begin{equation}\label{1.4}
 \lim_{l\to\infty}s\big(\mathcal J_{A, B_l(0)}\big)=s(A^\infty+D).\zzz
 \end{equation}
 \end{enumerate}
\end{introtheo}

In Theorem~\ref{th1.2}\textup{(ii)}, $A(x)$ need not be symmetric,
and convergence to $A^\infty$ is required only on the cubes $H_k$.
Examples~\ref{ex4.1} and~\ref{ex4.2} below show that the global upper bound \eqref{1.2} cannot in general be omitted and that convergence along a sequence of points does not suffice.
Example~\ref{ex4.4}, based on a translated Gaussian kernel, shows that the evenness assumption on the kernels cannot in general be omitted either.

If $A(0)$ is irreducible, Theorem~\ref{th1.2}\textup{(i)} implies $\lim_{l\searrow0}s\big(\mathcal J_{A, B_l(0)}\big)=s(A(0))$. Thus, if $s(A(0))<0<\lim_{l\to\infty}s\big(\mathcal J_{A, B_l(0)}\big)$, continuity implies $s\big(\mathcal J_{A, B_L(0)}\big)=0$ for some $L>0$. If zero is the principal eigenvalue for this domain, strict domain comparison implies that $L$ is unique.\vspace{-1mm}

\begin{introcoro}[A unique critical domain size]\label{cor1.1}
 Suppose that $A(0)$ is irreducible and that, for some $L>0$, zero is the principal eigenvalue of $\mathcal J_{A, B_L(0)}$. Then we have\zzz
  \begin{equation}\label{1.5}
 s\big(\mathcal J_{A, B_l(0)}\big)<0\;\;\;\text{for}\;\;0<l<L,\;\;\;
 s\big(\mathcal J_{A, B_l(0)}\big)>0\;\;\;\text{for}\;\;l>L.\zzz
    \end{equation}
Under the hypotheses of Theorem~\ref{th1.2}\textup{(ii)}, we also have $s(A^\infty+D)>0$.\zzz
\end{introcoro}

Example~\ref{ex4.3} below has symmetric kernels and a symmetric irreducible coefficient matrix, but its spectral bound is zero on every sufficiently large interval centered at the origin. Thus symmetry alone does not imply strict domain monotonicity or uniqueness of the domain size at which the spectral bound vanishes.

To address Question~(2), we consider two nonlinear models on one-dimensional habitats. Under the respective assumptions, both models admit a finite critical habitat length, defined as the infimum of the interval lengths for which the spectral bound is positive. Whenever the set of lengths with zero spectral bound is nonempty, the critical length is its right endpoint. The resulting dynamical thresholds remain valid even when the spectral bound vanishes on an interval of lengths. For the heterogeneous epidemic model, which is related to the reaction-diffusion model of Li, Li and Wang~\cite{LLW25SAM}, no spatial monotonicity of the coefficients is required. Under strong subhomogeneity and the existence of a positive upper solution, the results of Wang and Zhang~\cite{WZ26JMPA} ensure extinction at and below the critical length and a unique positive equilibrium above it, globally asymptotically stable for nontrivial nonnegative initial data. For the scalar logistic equation with a free boundary and a sign-changing growth rate, spreading occurs for every positive boundary parameter when the initial habitat length is at least the critical length. For smaller initial habitats, there is a positive finite threshold in the boundary parameter: vanishing occurs at and below the threshold, while spreading occurs above it. At the critical parameter, the limiting boundary position equals the critical habitat length. A zero spectral bound precludes exponential decay of the linearized semigroup in operator norm (see more details in Corollary\ref{cor5.1} and Theorem \ref{th6.1}). 

The paper is organized as follows. Section~\ref{sec2} addresses spectral characterizations and the existence criterion for principal eigenvalues. 
Section~\ref{sec3} establishes continuity with respect to coefficients and domains and proves Theorem~\ref{th1.1}. Section~\ref{sec4} proves Theorem~\ref{th1.2} and Corollary~\ref{cor1.1} and presents the counterexamples. Section~\ref{sec5} treats the epidemic model, and
Section~\ref{sec6} studies the scalar free-boundary problem.

\section{Spectral characterizations and principal eigenfunctions}\label{sec2}

For a bounded operator or a matrix $L$, denote its spectral radius by $r(L)=\sup\{|\lambda|:\lambda\in\sigma(L)\}$. For coefficient matrices defined only on $\boo$, we define $\mathcal J_{B,\oo}$ by \eqref{1.1} with $A$ replaced by $B$,
and define $s(\mathcal J_{B,\oo})$ and $s(B,\overline{\Omega})$ in the same way.
Matrix inequalities are entrywise: $A\le B$ means $a_{ij}(x)\le b_{ij}(x)$ on $\boo$ for all $i,j\in\mathcal S$. We write $A\le B$, $A\not\equiv B$ if, in addition, at least one
of these inequalities is strict at some point of $\boo$. We also use $\mathbb R_+=[0,\infty)$ and $t_+=\max\{t,0\}$ for $t\in\mathbb R$.

\subsection{Spectral characterizations}

The following theorem collects important spectral characterizations about $\mathcal J_{A,\Omega}$.\vspace{-1mm}

\begin{theo}[Spectral characterizations]\label{th2.1} The following statements are valid:\vspace{-2mm}
    \begin{enumerate}[label=\textup{(\roman*)},leftmargin=8mm]
 \item For every sufficiently small $\ep>0$, there exist continuous cooperative control matrices $A_\ep^-\le A\le A_\ep^+$ on $\boo$, with $A_\ep^\pm\to A$ uniformly on $\boo$ as $\ep\searrow0$. Each $\mathcal J_{A_\ep^\pm,\oo}$ has an eigenfunction $\varphi_\ep^\pm\gg 0$, and its principal eigenvalue $\lambda_p(\mathcal J_{A_\ep^\pm,\oo})$ satisfies\vspace{-1mm}
  \begin{equation}\label{2.1}
  \lambda_p\big(\mathcal J_{A_\ep^-,\oo}\big)\le s\big(\mathcal J_{A,\oo}\big)\le\lambda_p\big(\mathcal J_{A_\ep^+,\oo}\big),\;\;\;
  \lim_{\ep\searrow0}\lambda_p\big(\mathcal J_{A_\ep^\pm,\oo}\big)=s\big(\mathcal J_{A,\oo}\big).\vspace{-2mm}
  \end{equation} 
 \item The Collatz--Wielandt formulas are
 \begin{equation}\label{2.2}
  s(\mathcal J_{A,\oo})=\sup_{\varphi\gg0}\,\inf_{x\in\boo,\,i\in\mathcal S}
  \frac{\mathcal J_{i, A,\oo}[\varphi](x)}{\varphi_i(x)}
  =\inf_{\varphi\gg0}\,\sup_{x\in\boo,\,i\in\mathcal S}
  \frac{\mathcal J_{i, A,\oo}[\varphi](x)}{\varphi_i(x)}.
  \end{equation}
 Equivalently,
 \bee
 \begin{aligned}
 s(\mathcal J_{A,\oo})&=\sup\{\mu\in\mathbb R:\exists\,\varphi\gg0,\,
 \mathcal J_{A,\oo}[\varphi]\ge\mu\varphi\}\\
  &=\inf\{\mu\in\mathbb R:\exists\,\varphi\gg0,\,
   \mathcal J_{A,\oo}[\varphi]\le\mu\varphi\},
    \label{2.3}
  \end{aligned}
 \eee
 where the inequalities hold componentwise on $\boo$.
  \item One has
  \begin{equation}\label{2.4}
  s(A,\overline{\Omega})\le s(\mathcal J_{A,\oo}). \end{equation}
 Moreover, $s(\mathcal J_{A,\oo})$ is the principal eigenvalue if and only if $s(\mathcal J_{A,\oo})>s(A,\overline{\Omega})$.
 \item For $\varphi,\psi\gg0$, set
  \[\mathcal Q(\varphi,\psi)=\frac{\int_\oo\psi(x)^{\mathsf T}
  \mathcal J_{A,\oo}[\varphi](x)\dx}
 {\int_\oo\psi(x)^{\mathsf T}\varphi(x)\dx}.\]
 Then, without any symmetry assumption, we have
 \begin{equation}\label{2.5}
  s(\mathcal J_{A,\oo})=\sup_{\varphi\gg0}\,\inf_{\psi\gg0}
   \mathcal Q(\varphi,\psi)
   =\inf_{\varphi\gg0}\,\sup_{\psi\gg0}
   \mathcal Q(\varphi,\psi).\vspace{-2mm}
  \end{equation}   
  \item For a constant positive diagonal matrix $W=\diag(w_1,\ldots,w_m)$, define
  \[
 \mathcal R_W(\varphi)=\frac{\int_\oo\varphi(x)^{\mathsf T}
 W\mathcal J_{A,\oo}[\varphi](x)\dx}
{\int_\oo\varphi(x)^{\mathsf T}W\varphi(x)\dx},
 \;\;\; 0\ne\varphi\in[L^2(\oo)]^m.
  \]
Here all functions are real-valued, and $\mathcal J_{A,\oo}$ is defined on $[L^2(\oo)]^m$ by \eqref{1.1}. The notation $s(\mathcal J_{A,\oo})$ still refers to the spectral bound on $[C(\boo)]^m$. Then
  \bee s(\mathcal J_{A,\oo}) \le\sup_{0\ne\varphi\in[L^2(\oo)]^m}\mathcal R_W(\varphi).
 \label{2.6}\eee
Equality holds in \eqref{2.6} if $J_i(x,y)=J_i(y,x)$ for every component and $WA(x)=A(x)^{\mathsf T}W$ on $\boo$. In particular, taking $W=I$, we obtain the ordinary Rayleigh quotient $\mathcal R=\mathcal R_I$ and the corresponding formula when both the kernels and $A$ are symmetric.\vspace{-1mm}
\end{enumerate} 
\end{theo}

\begin{proof} \textup{(i)} This follows from \cite[Theorem A]{WZ26JMPA}.
The control matrices can be chosen to satisfy \vspace{-1mm}
\[A-2\ep I\le A_\ep^-\le A-\ep I,\;\;\;
A+\ep I\le A_\ep^+\le A+2\ep I,\;\;\; A_\ep^+=A_\ep^-+3\ep I.\vspace{-1mm}\]
They differ from $A$ only by scalar diagonal perturbations, and
$\lambda_p(\mathcal J_{A_\ep^+,\oo})=\lambda_p(\mathcal J_{A_\ep^-,\oo})+3\ep$.

\textup{(ii)} The Collatz--Wielandt formulas \eqref{2.2} follow from \cite[Theorem A]{WZ26JMPA}. 
To prove \eqref{2.3}, fix $\varphi\gg0$. Since every $\varphi_i$ has a positive
minimum on $\boo$, the ratios $\mathcal J_{i, A,\oo}[\varphi]/\varphi_i$ are continuous
there and attain their minima and maxima. Componentwise division and \eqref{2.2} give
  \begin{align*}
\mathcal J_{A,\oo}[\varphi]\ge\mu\varphi,
&\;\;\;\Longrightarrow\;
\mu\le\min_{x\in\boo,\,i\in\mathcal S}
\frac{\mathcal J_{i, A,\oo}[\varphi](x)}{\varphi_i(x)}\le s(\mathcal J_{A,\oo}),\\
\mathcal J_{A,\oo}[\varphi]\le\mu\varphi,
&\;\;\;\Longrightarrow\;
\mu\ge\max_{x\in\boo,\,i\in\mathcal S}
\frac{\mathcal J_{i, A,\oo}[\varphi](x)}{\varphi_i(x)}\ge s(\mathcal J_{A,\oo}).
\end{align*}
It follows that the supremum in \eqref{2.3} is at most $s(\mathcal J_{A,\oo})$,
and the infimum is at least $s(\mathcal J_{A,\oo})$.

To prove the reverse inequalities, we use the eigenfunctions in \textup{(i)}.
Since $A_\ep^-\le A\le A_\ep^+$,
\[\mathcal J_{A,\oo}[\varphi_\ep^-]\ge\mathcal J_{A_\ep^-,\oo}[\varphi_\ep^-]=
\lambda_p\big(\mathcal J_{A_\ep^-,\oo}\big)\varphi_\ep^-,\;\;\;
\mathcal J_{A,\oo}[\varphi_\ep^+]\le \mathcal J_{A_\ep^+,\oo}[\varphi_\ep^+]
=\lambda_p\big(\mathcal J_{A_\ep^+,\oo}\big)\varphi_\ep^+.\]
Hence $\lambda_p(\mathcal J_{A_\ep^-,\oo})$ belongs to the first set in
\eqref{2.3}, and $\lambda_p(\mathcal J_{A_\ep^+,\oo})$ belongs to the second. Letting $\ep\searrow0$, we deduce the reverse inequalities from \eqref{2.1}. This proves \eqref{2.3}.

\textup{(iii)} This follows from \cite[Lemma A.1(i)--(ii)]{WZ26JMPA}.

\textup{(iv)} \textit{Step 1: Extrema of $\mathcal Q(\varphi,\cdot)$.} Fix $\varphi\gg0$ and set
\[
q_i(x)=\frac{\mathcal J_{i, A,\oo}[\varphi](x)}{\varphi_i(x)},\;\;\;
\alpha=\min_{x\in\boo,\,i\in\mathcal S}q_i(x),\;\;\;
\beta=\max_{x\in\boo,\,i\in\mathcal S}q_i(x).
\]
The functions $q_i$ are continuous on $\boo$. We first prove that
  \begin{equation}\label{2.7}
  \inf_{\psi\gg0}\mathcal Q(\varphi,\psi)=\alpha,\;\;\;
  \sup_{\psi\gg0}\mathcal Q(\varphi,\psi)=\beta.\end{equation}
  
For every $\psi\gg0$, the weights $\psi_i\varphi_i$ are positive, so
 \[\alpha\le\mathcal Q(\varphi,\psi)=\frac{\sum_i\int_\oo\psi_i\varphi_i q_i\dx}
   {\sum_i\int_\oo\psi_i\varphi_i\dx}\le\beta.\]
To prove the first equality in \eqref{2.7}, choose $x_*\in\boo$ and $i_0\in\mathcal S$ such that $q_{i_0}(x_*)=\alpha$. Given $\eta>0$, by continuity there exists $r>0$ such that
  \[q_{i_0}(x)\le\alpha+\eta\;\;\;\text{for all }x\in U:=B_r(x_*)\cap\oo.\]
Since $x_*\in\boo$, the set $U$ is nonempty and open, even if $x_*\in\partial\oo$.
Choose a nonzero, nonnegative function $h\in C_c(U)$ and extend it
by zero to $\boo$. Since $\varphi_{i_0}>0$, we may rescale $h$ so that $\int_\oo h\varphi_{i_0}\dx=1$. For $\delta>0$, set $\psi^\delta=h\mathbf e_{i_0}+\delta\mathbf1$, 
where $\mathbf e_{i_0}$ is the $i_0$-th  standard basis vector and
$\mathbf1=(1,\ldots,1)^{\mathsf T}$. Each component of $\psi^\delta$
is continuous on $\boo$ and at least $\delta$, so $\psi^\delta\gg0$.
Substituting this function into $\mathcal Q$, we obtain
\begin{align*}
\mathcal Q(\varphi,\psi^\delta)
&=\frac{\int_\oo h\varphi_{i_0}q_{i_0}\dx+\delta\sum_i\int_\oo\varphi_iq_i\dx}
{1+\delta\sum_i\int_\oo\varphi_i\dx}
\to\int_\oo h\varphi_{i_0}q_{i_0}\dx
\le\alpha+\eta\;\;\;(\delta\searrow0).
\end{align*}
Consequently,
 \[\alpha\le\inf_{\psi\gg0}\mathcal Q(\varphi,\psi)
\le\lim_{\delta\searrow0}\mathcal Q(\varphi,\psi^\delta)\le\alpha+\eta.\]
Letting $\eta\searrow0$ proves the first equality.

To prove the second equality in \eqref{2.7}, fix $\eta>0$. By continuity at a maximum point, choose an index $j_0$ and a nonempty open set $V\subset\oo$ such that
$q_{j_0}\ge\beta-\eta$ on $V$. Choose $g\in C_c(V)$ with $g\ge0$ and $\int_\oo g\varphi_{j_0}\dx=1$, and extend $g$ by zero to $\boo$. For $\delta>0$, set $\widetilde\psi^\delta=g\mathbf e_{j_0}+\delta\mathbf1\gg0$. The same substitution gives
\begin{align*}
\beta-\eta
&\le\int_\oo g\varphi_{j_0}q_{j_0}\dx
=\lim_{\delta\searrow0}
 \mathcal Q(\varphi,\widetilde\psi^\delta)
\le\sup_{\psi\gg0}\mathcal Q(\varphi,\psi)\le\beta.
\end{align*}
Letting $\eta\searrow0$ proves the second equality.

\textit{Step 2: Dual Collatz--Wielandt formulas.} Since $\varphi\gg0$ was arbitrary, the first equality in \eqref{2.7} and the sup-inf formula in \eqref{2.2} yield
  \[\sup_{\varphi\gg0}\,\inf_{\psi\gg0}\mathcal Q(\varphi,\psi)
=\sup_{\varphi\gg0}\,\min_{x\in\boo,\,i\in\mathcal S}
 \frac{\mathcal J_{i, A,\oo}[\varphi](x)}{\varphi_i(x)}
=s(\mathcal J_{A,\oo}).\]
The second equality and the inf-sup formula in \eqref{2.2} yield
  \[\inf_{\varphi\gg0}\,\sup_{\psi\gg0}\mathcal Q(\varphi,\psi)
=\inf_{\varphi\gg0}\,\max_{x\in\boo,\,i\in\mathcal S}
 \frac{\mathcal J_{i, A,\oo}[\varphi](x)}{\varphi_i(x)}=s(\mathcal J_{A,\oo}).\]
This proves \eqref{2.5}.

\textup{(v)} \textit{Step 1: Weighted Rayleigh bound.} Write 
 \[\langle u,v\rangle_W=\int_\oo u^{\mathsf T}Wv\dx,\;\;\;
\|v\|_W^2=\langle v,v\rangle_W.\]
For the weighted Rayleigh quotient in part~\textup{(v)}, the supremum
is unchanged when restricted to continuous or strongly positive
continuous functions:\vspace{-1mm}
\begin{align}
    \sup_{\varphi\gg0}\mathcal R_W(\varphi)
    =\sup_{0\ne\varphi\in [C(\boo)]^m}\mathcal R_W(\varphi)
    =\sup_{0\ne\varphi\in[L^2(\oo)]^m}\mathcal R_W(\varphi).
    \label{2.8}
\end{align}
Indeed, the last two suprema are equal by the density of $[C(\boo)]^m$
in $[L^2(\oo)]^m$ and the continuity of $\mathcal R_W$ away from zero.
For $0\ne v\in[C(\boo)]^m$, the nonnegativity of the kernels and the
cooperativity of $A$ imply $\mathcal R_W(v)\le\mathcal R_W(|v|)$, where
$|v|$ is taken componentwise. Since $|v|+\delta\mathbf1\gg0$ and $\mathcal R_W(|v|+\delta\mathbf1)\to\mathcal R_W(|v|)$ as $\delta\searrow0$, the first equality follows as well.

Since $\mathcal J_{A,\oo}[\varphi_\ep^-]\ge\lambda_p\big(\mathcal J_{A_\ep^-,\oo}\big)\varphi_\ep^-$, we have
  \[\mathcal R_W(\varphi_\ep^-)=\frac{\int_\oo\varphi_\ep^-(x)^{\mathsf T}W\mathcal J_{A,\oo}[\varphi_\ep^-](x)\dx} {\int_\oo\varphi_\ep^-(x)^{\mathsf T}W\varphi_\ep^-(x)\dx}
  \ge\lambda_p\big(\mathcal J_{A_\ep^-,\oo}\big)\to s(\mathcal J_{A,\oo}).\]
Together with \eqref{2.8}, this proves \eqref{2.6}. 

\textit{Step 2: Equality in \eqref{2.6}.} Suppose that $J_i(x,y)=J_i(y,x)$ and $WA=A^{\mathsf T}W$. Let $u=\varphi_\ep^+$. Then $u\gg0$ and $\mathcal J_{A_\ep^+,\oo}[u]=\lambda_p\big(\mathcal J_{A_\ep^+,\oo}\big)u$. For any $v\in[L^2(\oo)]^m$, we have
 \bee
2v_i(x)v_j(y)\le\frac{u_j(y)}{u_i(x)}v_i^2(x)+\frac{u_i(x)}{u_j(y)}v_j^2(y).\label{2.9}\eee
Taking $j=i$ in \eqref{2.9}, integrating against $J_i(x,y)$ over $\Omega\times\Omega$, and using $J_i(x,y)=J_i(y,x)$, we obtain
\begin{equation}\label{2.10}
\int_\Omega\int_\Omega J_i(x,y)v_i(x)v_i(y)\dy\dx
\le\int_\Omega\frac{v_i^2(x)}{u_i(x)}
\int_\Omega J_i(x,y)u_i(y)\dy\dx.
\end{equation}
Taking $y=x$ in \eqref{2.9}, summing over $i\ne j$ with the nonnegative weights $w_i a_{ij}(x)$, and using $w_i a_{ij}=w_j a_{ji}$, we obtain
 \begin{align}
\sum_{i=1}^m\sum_{j\ne i}w_i a_{ij}(x)v_i(x)v_j(x)&\le\frac12\sum_{i=1}^m\sum_{j\ne i}w_i a_{ij}(x)
 \left(\frac{u_j(x)}{u_i(x)}v_i^2(x)+\frac{u_i(x)}{u_j(x)}v_j^2(x)\right)\nonumber\\
&=\frac12\sum_{i=1}^m\sum_{j\ne i} \bigl(w_i a_{ij}(x)+w_j a_{ji}(x)\bigr)
 \frac{u_j(x)}{u_i(x)}v_i^2(x)\nonumber\\
&=\sum_{i=1}^m\frac{w_i v_i^2(x)}{u_i(x)}\sum_{j\ne i}a_{ij}(x)u_j(x).
\label{2.11}
\end{align}
Note that $a_{ij}(x)=a_{\ep,ij}^+(x)$ when $j\ne i$, and $a_{ii}(x)\le a_{\ep,ii}^+(x)$. Integrating \eqref{2.11} over $\oo$, adding the diagonal terms, and using \eqref{2.10} multiplied by $w_i d_i$ and summed over $i$,  we obtain
\begin{align*}
\mathcal R_W(v)\|v\|_W^2&=\big\langle v,\mathcal J_{A,\oo}[v]\big\rangle_W\\
&\le\sum_{i=1}^m\int_\oo\frac{w_i v_i^2(x)}{u_i(x)}
 \biggl[d_i\int_\oo J_i(x,y)u_i(y)\dy+ a^+_{\ep,ii}(x)u_i(x)
 +\sum_{j\ne i}a_{ij}(x)u_j(x)\biggr]\dx\\
&=\sum_{i=1}^m\int_\oo
 w_i v_i^2(x)\frac{\mathcal J_{i, A_\ep^+,\oo}[u](x)}{u_i(x)}\dx
=\lambda_p\big(\mathcal J_{A_\ep^+,\oo}\big)\|v\|_W^2.
\end{align*}
Hence $\mathcal R_W(v)\le\lambda_p\big(\mathcal J_{A_\ep^+,\oo}\big)$ for all $0\ne v\in[L^2(\oo)]^m$. Letting $\varepsilon\searrow0$ and using \eqref{2.1}, we obtain
 the reverse inequality in \eqref{2.6}.
\www\end{proof}

The weighted version of \eqref{2.10} uses the constancy of $W$. For a variable weight
$W(x)$, symmetry of $J_i$ alone does not imply
$w_i(x)J_i(x,y)=w_i(y)J_i(y,x)$, which is needed to interchange $x$
and $y$ in the weighted integral.

\subsection{A sufficient condition for a principal eigenvalue}

For scalar operators, Coville~\cite[Theorem 1.1]{C10JDE} established a nonintegrability condition for the existence of a principal eigenvalue. For cooperative systems, Wang and Zhang~\cite[Lemma 2.1(ii)]{WZ26JMPA} proved a related criterion assuming a strongly positive eigenvector at a
maximizer of the spectral bound of the coefficient matrix. The following proposition requires irreducibility at only one point of $\boo$, which need not be a maximizing point.\vspace{-1mm}

\begin{prop}[Principal eigenfunctions]\label{prop2.1} Assume \textbf{\bf(J)} and \textbf{\bf(A)} hold. If
\begin{equation}\label{2.12}
\int_\oo\frac{\dx}{s(A,\overline{\Omega})-s(A(x))}=\infty,
\end{equation} 
where $1/0=\infty$, then $s(\mathcal J_{A,\oo})> s(A,\overline{\Omega})$ and
 $s(\mathcal J_{A,\oo})$ is the principal eigenvalue of $\mathcal J_{A,\oo}$. 
In particular, this conclusion holds if $s(A(x))\equiv s(A,\overline{\Omega})$ on a nonempty open
subset of $\oo$.\www
\end{prop}

\begin{proof} By \textbf{(J)} and the compactness of $\boo$, there exist $r,\kappa>0$ such that\zzz
    \[ J_i(x,y)\ge\kappa\;\;\;\text{whenever }x,y\in\boo,\ |x-y|<2r,
    \;\; i\in\mathcal S.\zzz\]
Cover $\boo$ by finitely many balls of radius $r$ with centers in $\boo$. It follows from \eqref{2.12} that at least one of the sets $U=\oo\cap B_r(z)$ for some $z\in\boo$ satisfies
    \bee
    \int_U\frac{\dx}{s(A,\overline{\Omega})-s(A(x))}=\infty.
    \label{2.13}
    \eee
        
In view of \eqref{2.4}, it remains to exclude $s(\mathcal J_{A,\oo})=s(A,\overline{\Omega})$. Suppose, to the contrary, that equality holds. For $\ep>0$, set
    \[ R^\ep(x)=\kk((s(A,\overline{\Omega})+\ep)I-A(x)\rr)^{-1}\ge0.\]
The inverse is nonnegative since $A(x)$ is cooperative and $s(A,\overline{\Omega})+\ep>s(A(x))$.
    
\textit{Claim.} For every $\ep>0$,
    \begin{equation}\label{2.14}
        \int_U\operatorname{tr}R^\ep(x)\dx
        \le\frac1\kappa\sum_{i=1}^m\frac1{d_i}.
    \end{equation}
Fix $\varepsilon>0$. By the inf-sup characterization \eqref{2.2}, there exists $\varphi^\varepsilon\gg0$ satisfying $\mathcal J_{A,\Omega}[\varphi^\varepsilon]\le \big(s(A,\overline{\Omega})+\varepsilon\big)\varphi^\varepsilon$. 
Multiply this inequality by $R^\varepsilon(x)$. Since all entries of $R^\varepsilon$ are nonnegative, we obtain, for every $x\in U$ and each index $i$,
  \[\varphi_i^\varepsilon(x)\ge\sum_{j=1}^m R^\varepsilon_{ij}(x)d_j\int_\Omega J_j(x,y)\varphi_j^\varepsilon(y)\dy
  \ge\kappa d_i R^\varepsilon_{ii}(x)\int_U\varphi_i^\varepsilon(y)\dy.\]
Integration over $U$ and division by $\int_U\varphi_i^\ep(y)\dy>0$ yield $\int_U R^\ep_{ii}(x)\dx\le\frac1{\kappa d_i}$. Summing over $i$ proves the claim.    

The spectral bound of a cooperative matrix is a real eigenvalue. For each fixed $x$, denote the eigenvalues of $A(x)$, counted with algebraic multiplicity, by $\lambda_1(x)=s(A(x)),\lambda_2(x),\ldots,\lambda_m(x)$. Thanks to $\operatorname{Re}\lambda_j(x)\le s(A(x))$, we have
 \[\operatorname{Re}\frac1{s(A,\overline{\Omega})+\ep-\lambda_j(x)}=
 \frac{s(A,\overline{\Omega})+\ep-\operatorname{Re}\lambda_j(x)}
 {|s(A,\overline{\Omega})+\ep-\lambda_j(x)|^2}>0.\]
Since the trace is the sum of the eigenvalues counted with algebraic multiplicity, we have
 \[\operatorname{tr}R^\ep(x)=\sum_{j=1}^m\operatorname{Re}
 \frac1{s(A,\overline{\Omega})+\ep-\lambda_j(x)}\ge\frac1{s(A,\overline{\Omega})+\ep-s(A(x))}.\]
Combining this estimate with \eqref{2.14}, we obtain
  \[\int_U\frac{\dx}{s(A,\overline{\Omega})+\ep-s(A(x))} \le\frac1\kappa\sum_{i=1}^m\frac1{d_i}.\]
Together with \qq{2.13}, this gives
 \[\infty=\int_U\frac{\dx}{s(A,\overline{\Omega})-s(A(x))}=\lim_{\ep\searrow0}
    \int_U\frac{\dx}{s(A,\overline{\Omega})+\ep-s(A(x))}
    \le\frac1\kappa\sum_{i=1}^m\frac1{d_i}<\infty,    \]
a contradiction. Hence $s(\mathcal J_{A,\oo})>s(A,\overline{\Omega})$, and the existence of the principal eigenvalue follows from Theorem~\ref{th2.1}\textup{(iii)}. The last assertion follows since the integrand in \eqref{2.12} is infinite on a set of positive measure.\zzz
\end{proof}

In one dimension, local Lipschitz continuity of $s(A(x))$ at a
maximizing point implies \eqref{2.12}.

\begin{coro}\label{cor2.1}
Assume {\bf(J)} and {\bf(A)} hold, and let $\oo=(a,b)$.
If $x\mapsto s(A(x))$ is locally Lipschitz near a point
$x_*\in[a,b]$ where it attains its maximum, with the neighborhood
taken relative to $[a,b]$, then $\mathcal J_{A,\oo}$  has a principal eigenvalue.\zzz
\end{coro}

\begin{proof} By local Lipschitz continuity, there exist $C,r>0$ such that
\[0\le  s(A,\overline{\Omega})-s(A(x))\le C|x-x_*|,\;\;\; x\in(a,b)\cap B_r(x_*).\]
Therefore
  \[\int_\oo\frac{\dx}{ s(A,\overline{\Omega})-s(A(x))}
  \ge\frac1C\int_{(a,b)\cap B_r(x_*)}\frac{\dx}{|x-x_*|}=\infty.\vspace{0.5mm}\]
The conclusion follows from Proposition \ref{prop2.1}.\zzz
\end{proof}

\section{Continuity and monotonicity}\label{sec3}

The following strong maximum principle, which follows from the arguments in \cite[Lemma 3.1 and Theorem 3.1]{WZ26JMPA}, will be used below. \zzz

\begin{lem}\label{lem3.1}
Assume \textbf{\bf(J)} holds.  Let $B=(b_{ij})\in[C(\overline {\oo})]^{m\times m}$ be cooperative, with $B(x_*)$ irreducible for some $x_*\in\overline\oo$. Suppose that $v\in[L^\infty(\oo)]^m$, $v\ge0$, and, for some $\lambda\in\mathbb R$,\zzz
\[
d_i\int_{\oo} J_i(x,y)v_i(y)\dy
+\sum_{j=1}^m b_{ij}(x)v_j(x)\le\lambda v_i(x)
\;\;\;\text{for a.e. }x\in \oo,\;\;\; i\in\mathcal S.\www
\]
Then either $v=0$ a.e. in $\oo$, or $\operatorname*{ess\,inf}_{x\in \oo}v_i(x)>0$ for every $i\in\mathcal S$. 
In particular, if $v\in [C(\overline {\oo})]^m$, then either $v\equiv0$ or $v\gg0$.
\end{lem}

\subsection{Dependence on the coefficients---Proof of Theorem~\ref{th1.1}\textup{(i)}}

We first prove monotonicity and continuity with respect to the coefficients. Neither result requires the existence of a principal eigenvalue.\zzz

\begin{lem}\label{lem3.2} Assume \textbf{\bf(J)} and \textbf{\bf(A)} hold. \www
\begin{enumerate}[label=\textup{(\roman*)},leftmargin=8mm]
\item Let $B$ satisfy \textbf{\bf(A)}, with an irreducibility point in
$\boo$. If $A\ge B$ on $\boo$, then $s(\mathcal J_{A,\oo})\ge s(\mathcal J_{B,\oo})$.\zzz
\item Let each $A_k$ satisfy \textbf{\bf(A)}, with an irreducibility point in $\boo$. If $A_k\to A$ uniformly on $\boo$, then $\lim_{k\to\infty}s(\mathcal J_{A_k,\oo})=s(\mathcal J_{A,\oo})$.\www
\end{enumerate}
\end{lem}\www

\begin{proof}  Part~\textup{(i)} follows from \qq{2.2}.

For part~\textup{(ii)}, fix a small $\ep>0$, and let $\varphi_\ep^\pm\gg 0$ be the eigenfunctions of the control operators associated with $A_\ep^\pm$ in Theorem~\ref{th2.1}. Then we have
  \[\delta_{k,\ep}:=\max_{\tau\in\{-,+\}}\max_{x\in\boo,\,i\in\mathcal S}
\frac{|((A_k(x)-A(x))\varphi_\ep^\tau(x))_i|}{\varphi_{\ep,i}^\tau(x)}\to0
\;\;\;(k\to\infty).\]
The inequalities $\mathcal J_{A,\oo}[\varphi_\ep^-]\ge\lambda_p\big(\mathcal J_{A_\ep^-,\oo}\big)\varphi_\ep^-$ and $\mathcal J_{A,\oo}[\varphi_\ep^+]\le\lambda_p\big(\mathcal J_{A_\ep^+,\oo}\big)\varphi_\ep^+$ imply\zzz
  \[\lambda_p\big(\mathcal J_{A_\ep^-,\oo}\big) -\delta_{k,\ep}
\le s(\mathcal J_{A_k,\oo})\le\lambda_p\big(\mathcal J_{A_\ep^+,\oo}\big) +\delta_{k,\ep}\]
by \eqref{2.2}. Letting $k\to\infty$ and then $\ep\searrow0$, we obtain
$s(\mathcal J_{A_k,\oo})\to s(\mathcal J_{A,\oo})$.\www
\end{proof}

We remark that the above lemma still true when the irreducibility property fails (see, e.g., \cite{Z24ANS}).
\begin{proof}[Proof of Theorem~\ref{th1.1}\textup{(i)}] 

\textit{Case 1: Both operators have principal eigenvalues}.  Choose $u,v\gg0$ with
$\mathcal J_{A,\oo}[u]=s(\mathcal J_{A,\oo})u$ and $\mathcal J_{B,\oo}[v]=s(\mathcal J_{B,\oo})v$.
By Lemma~\ref{lem3.2}\textup{(i)}, we have $s(\mathcal J_{B,\oo})\le s(\mathcal J_{A,\oo})$. Suppose, to the contrary, that equality holds. Set
  \[t=\max_{x\in\boo,\,i\in\mathcal S}\frac{v_i(x)}{u_i(x)},\;\;\; w=tu-v\ge0.\]
The maximum defining $t$ is attained at some $(\bar x,i_0)$, so $w_{i_0}(\bar x)=0$. The identity
  \[\mathcal J_{A,\oo}[w]-s(\mathcal J_{A,\oo})w=-(A-B)v\le0\]
and Lemma~\ref{lem3.1} imply $w\equiv0$. This contradicts $A\ge B$, $A\not\equiv B$, and $v\gg0$. Thus $s(\mathcal J_{A,\oo})>s(\mathcal J_{B,\oo})$.

\textit{Case 2: Only $\mathcal J_{A,\oo}$ has a principal eigenvalue}. By Theorem~\ref{th2.1}, we have
\[
s(\mathcal J_{B,\oo})=s(B,\overline{\Omega})\le s(A,\overline{\Omega})<s(\mathcal J_{A,\oo}).
\]

\textit{Case 3: Only $\mathcal J_{B,\oo}$ has a principal eigenvalue}. Set $A_t=B+t(A-B)$ for $0<t<1$. Then $B\le A_t\le A$, so
$A_t$ is cooperative and is irreducible at every irreducibility
point of $B$.
Since $A_t\to B$ uniformly as $t\searrow0$,
\[
s(A_t,\overline{\Omega})\to s(B,\overline{\Omega})<s(\mathcal J_{B,\oo}).
\]
Choose $t>0$ small enough that $s(A_t,\overline{\Omega})<s(\mathcal J_{B,\oo})\le s(\mathcal J_{A_t,\oo})$.
By Theorem~\ref{th2.1}, the operator $\mathcal J_{A_t,\oo}$ has a principal eigenvalue.  Case 1 with $B$ and $A_t$ now gives $s(\mathcal J_{B,\oo})<s(\mathcal J_{A_t,\oo})\le s(\mathcal J_{A,\oo})$. The proof is complete.\www\end{proof}

If neither operator has a principal eigenvalue,  Theorem~\ref{th2.1}\textup{(iii)} implies the following conclusion.

\begin{coro}\label{cor3.1}
Let $A$ and $B$ satisfy \textbf{\bf(A)} with irreducibility points in $\oo$, and assume $A\ge B$ and $A\not\equiv B$.
Suppose that neither $\mathcal J_{A,\oo}$ nor
$\mathcal J_{B,\oo}$ has a principal eigenvalue. Then\zzz
\[
s(\mathcal J_{A,\oo})=s(A,\overline{\Omega}),\;\;\; s(\mathcal J_{B,\oo})=s(B,\overline{\Omega}).\www
\]
Thus strict inequality between the spectral bounds is equivalent to strict inequality between these maxima.\end{coro}

\subsection{Dependence on the domain---Proof of Theorem~\ref{th1.1}\textup{(ii)}}

The following proposition presents conditions for the convergence of the spectral bound as the domains increase or decrease to a limiting domain.\zzz

\begin{prop}\label{prop3.1}  The spectral bound $s(\mathcal J_{A,\oo})$
is continuous with respect to $\oo$ in the following sense:\zzz
  \begin{enumerate}[leftmargin=7mm]
\item[\textup{(i)}] if $\oo_k\subset\oo_{k+1}\subset\oo$ and $\oo_k\nearrow\oo$, then\zzz 
\begin{equation}\label{3.1}
    \lim_{k\to\infty}s(\mathcal J_{A,\oo_k})=s(\mathcal J_{A,\oo});\www
   \end{equation}
\item[\textup{(ii)}] if $\oo_k\supset\oo_{k+1}\supset\oo$,
$|\oo_k\setminus\oo|\to0$, and there are continuous maps
$h_k:\overline\oo_k\to\overline\oo$ such that
$h_k|_{\overline\oo}=\mathrm{id}$ and
$\sup_{x\in\overline\oo_k}|h_k(x)-x|\to0$,
then \eqref{3.1} still holds.
 \end{enumerate}
\end{prop}

\begin{proof} Let $G_1\subset G_2$ be bounded domains,  with boundaries of Lebesgue measure zero and an irreducibility point of $A$ in $\overline G_1$. We first prove
\begin{equation}\label{3.2}
s(\mathcal J_{A,G_1})\le s(\mathcal J_{A,G_2}).
\end{equation}
For any $\varphi\in[C(\overline G_2)]^m$ with $\varphi\gg0$, set
$\psi=\varphi|_{\overline G_1}\in[C(\overline G_1)]^m$.
Then $\psi\gg0$, and for $x\in\overline G_1$ and $i\in\mathcal S$ we have
\[
\mathcal J_{i, A,G_1}[\psi](x)=\mathcal J_{i, A,G_2}[\varphi](x)
-d_i\int_{G_2\setminus G_1}\!\!J_i(x,y)\varphi_i(y)\dy
\le\mathcal J_{i, A,G_2}[\varphi](x).\]
Since $\psi_i=\varphi_i>0$ on $\overline G_1$, the inf-sup formula \eqref{2.2} gives
\[s(\mathcal J_{A,G_1})\le\max_{\substack{x\in\overline G_1,\,i\in\mathcal S}}
 \frac{\mathcal J_{i, A,G_1}[\psi](x)}{\psi_i(x)}
\le\max_{\substack{x\in\overline G_2,\,i\in\mathcal S}}
 \frac{\mathcal J_{i, A,G_2}[\varphi](x)}{\varphi_i(x)}.\]
Taking the infimum over $\varphi$ proves \eqref{3.2} by \eqref{2.2} on $G_2$.

\textup{(i)} If $s(\mathcal J_{A,\oo})$ is not the principal eigenvalue, then $s(\mathcal J_{A,\oo})=s(A,\overline{\Omega})$. From \eqref{3.2}, we have $s(\mathcal J_{A,\oo_k})\le s(\mathcal J_{A,\oo})$, and Theorem~\ref{th2.1} yields 
  \[s(A,\overline{\Omega}_k)\le s(\mathcal J_{A,\oo_k})\le s(\mathcal J_{A,\oo}).\]
For any $\eta>0$, by continuity there exists $x_\eta\in\oo$ such that
$s(A(x_\eta))>s(A,\overline{\Omega})-\eta$. Since $\bigcup_k\oo_k=\oo$, we have
$x_\eta\in\oo_k$ for all sufficiently large $k$. Hence
$s(A,\overline{\Omega})-\eta<s(A,\overline{\Omega}_k)\le s(A,\overline{\Omega})$ for  these  $k$. Letting $\eta\searrow0$, we obtain $s(A,\overline{\Omega}_k)\to s(A,\overline{\Omega})$, and hence $s(\mathcal J_{A,\oo_k})\to s(\mathcal J_{A,\oo})$.

If $s(\mathcal J_{A,\oo})$ is the principal eigenvalue, then $s(\mathcal J_{A,\oo})>s(A,\overline{\Omega})$, and there exists $\varphi\gg0$ such that
$\mathcal J_{A,\oo}[\varphi]=s(\mathcal J_{A,\oo})\varphi$. The restriction of $\varphi$ to
$\overline\oo_k$ is continuous and strongly positive. For $x\in\overline\oo_k$,
 \begin{align}
\mathcal J_{i, A,\oo_k}[\varphi](x)=s(\mathcal J_{A,\oo})\varphi_i(x)
-d_i\int_{\oo\setminus\oo_k}J_i(x,y)\varphi_i(y)\dy\ge\bigl(s(\mathcal J_{A,\oo})-C|\oo\setminus\oo_k|\bigr)\varphi_i(x),
\label{3.3}\end{align}
where
  \[C=\max_{i\in\mathcal S}\frac{d_i\|J_i\|_\infty\|\varphi_i\|_{L^\infty(\oo)}}
  {\min_{\boo}\varphi_i}<\infty.\]
The sup-inf formula \eqref{2.2}, combined with \eqref{3.3}, implies $s(\mathcal J_{A,\oo_k})\ge s(\mathcal J_{A,\oo})-C|\oo\setminus\oo_k|$. Since $|\oo\setminus\oo_k|\to0$, the conclusion follows.

\textup{(ii)} For any $\ep>0$, by \eqref{2.2} there exists $\varphi^\ep\gg0$ such that
  \[\mathcal J_{A,\oo}[\varphi^\ep]\le(s(\mathcal J_{A,\oo})+\ep)\varphi^\ep
\;\;\;\text{on }\boo.\]
Set $\widehat\varphi_k^\ep=\varphi^\ep\circ h_k\in[C(\overline\oo_k)]^m$.
Then $\widehat\varphi_k^\ep\gg0$.
Let $E_k=\oo_k\setminus\oo$. Since $h_k|_{\boo}=\mathrm{id}$,
for $x\in\boo_k$ we have
 \begin{align}
\mathcal J_{i, A,\oo_k}[\widehat\varphi_k^\ep](x)
&=\mathcal J_{i, A,\oo}[\varphi^\ep](h_k(x))+d_i\int_\oo
 [J_i(x,y)-J_i(h_k(x),y)]\varphi_i^\ep(y)\dy\nonumber\\
&\;\;\;+[A_{i,*}(x)-A_{i,*}(h_k(x))]\varphi^\ep(h_k(x))+d_i\int_{E_k}J_i(x,y)
  \varphi_i^\ep(h_k(y))\dy.\label{3.4}
\end{align}
Each $J_i$ is uniformly continuous on $\boo_1\times\boo_1$,
and each $a_{ij}$ is uniformly continuous on $\boo_1$.
Since $\sup_{x\in\boo_k}|h_k(x)-x|\to0$, we obtain
\begin{align*}
\max_i\sup_{x\in\boo_k,\,y\in\boo}|J_i(x,y)-J_i(h_k(x),y)|\to0,\;\;\;
\max_{i,j}\,\sup_{x\in\boo_k}|a_{ij}(x)-a_{ij}(h_k(x))|\to0.
\end{align*}
Note that $\varphi^\ep$ is fixed and bounded. The second and third terms on the right-hand side of \eqref{3.4} tend to zero uniformly. The absolute value of the last term is at most
$d_i\|J_i\|_\infty\|\varphi_i^\ep\|_\infty|E_k|\to0$. Put $c_\ep=\min_{x\in\boo,\,i\in\mathcal S}\varphi_i^\ep(x)>0$. For all sufficiently large $k$, the sum of these three terms has
absolute value at most $\ep c_\ep\le\ep\widehat\varphi_{k,i}^\ep(x)$. Consequently,
  \[\mathcal J_{i,A,\Omega_k}[\widehat\varphi_k^\varepsilon](x) \le (s(\mathcal J_{A,\oo})+\varepsilon)\widehat\varphi_{k,i}^\varepsilon(x) +\varepsilon c_\varepsilon \le (s(\mathcal J_{A,\oo})+2\varepsilon)\widehat\varphi_{k,i}^\varepsilon(x). \]
By \eqref{2.2} and \eqref{3.2},  $s(\mathcal J_{A,\oo})\le s(\mathcal J_{A,\oo_k})\le s(\mathcal J_{A,\oo})+2\ep$ for all sufficiently large $k$. Letting $k\to\infty$ and then $\ep\searrow0$ proves \eqref{3.1}.\www
\end{proof}

In particular, Proposition \ref{prop3.1} applies to nested intervals, concentric balls, and centered cubes.

\begin{proof}[Proof of Theorem~\ref{th1.1}\textup{(ii)}] Since $|\oo_2\setminus\oo_1|>0$ and $J_i(x,x)>0$, it is easy to see that 
\begin{equation}\label{3.5}
\int_{\oo_1}\int_{\oo_2\setminus\oo_1}J_i(x,y)\dy\dx>0,
\;\;\;\forall\; i\in\mathcal S.
\end{equation}

\textit{Case 1: Principal eigenvalues on both domains}. 
Suppose both operators admit principal eigenvalues $\lambda_p(\mathcal J_{A,\oo_1})$ and $\lambda_p(\mathcal J_{A,\oo_2})$. Then $\lambda_p(\mathcal J_{A,\oo_1})\le \lambda_p(\mathcal J_{A,\oo_2})$ by \eqref{3.2}. Suppose, to the contrary, that $\lambda_p(\mathcal J_{A,\oo_1})=\lambda_p(\mathcal J_{A,\oo_2})=\lambda$. Let $u\in[C(\overline{\Omega}_1)]^m$ and $v\in[C(\overline{\Omega}_2)]^m$ be the corresponding eigenfunctions satisfying $u,v\gg0$, and set\vspace{-1mm}
 \[
t=\max_{x\in\overline\oo_1,\,i\in\mathcal S}
\frac{u_i(x)}{v_i(x)},\;\;\;
w=t v|_{\overline\oo_1}-u\ge0. \vspace{-1mm}
\]
At least one component of $w$ vanishes at a point of $\overline\oo_1$. For $x\in\overline\oo_1$,
\begin{equation}\label{3.6}
\mathcal J_{i, A,\oo_1}[w](x)-\lambda w_i(x)=-td_i\int_{\oo_2\setminus\oo_1}J_i(x,y)v_i(y)\dy\le0.
\end{equation}
Since $A$ has an irreducibility point in $\oo_1$,
Lemma~\ref{lem3.1} implies $w\equiv0$. It follows from \eqref{3.6} that $\int_{\oo_2\setminus\oo_1}J_i(x,y)v_i(y)\dy=0$ for all $x\in\oo_1$ and all $i$.
Since $\min_{\overline\oo_2}v_i>0$, integration in $x$ contradicts
\eqref{3.5}. Hence $\lambda_p(\mathcal J_{A,\oo_1})<\lambda_p(\mathcal J_{A,\oo_2})$.

\textit{Case 2: A principal eigenvalue only on $\oo_2$}. By Theorem~\ref{th2.1}, we have
  \[s(\mathcal J_{A,\oo_1})=s(A,\overline\Omega_1)\le s(A,\overline\Omega_2)<\lambda_p(\mathcal J_{A,\oo_2}).\]

\textit{Case 3: A principal eigenvalue only on $\oo_1$}. 
We have $s(A,\overline\Omega_1)<\lambda_p(\mathcal J_{A,\oo_1})$, and there exists $\eta>0$ small enough that
  \[s(A(x))<\lambda_p(\mathcal J_{A,\oo_1})\;\;\;\text{for all }x\in\overline\oo_2\;
\text{ with }\operatorname{dist}(x,\oo_1)\le\eta.\]
Let $\oo_\eta$ be the connected component containing $\oo_1$ of
$\oo_2\cap\{x:\operatorname{dist}(x,\oo_1)<\eta\}$. Then $|\oo_\eta\setminus\oo_1|>0$ when $\eta>0$ is sufficiently small. Thus, by \eqref{3.2}, $s(A,\overline\Omega_\eta)<\lambda_p(\mathcal J_{A,\oo_1})\le s(\mathcal J_{A,\oo_\eta})$. By Theorem~\ref{th2.1}(iii), a principal eigenvalue exists on $\oo_\eta$. Case 1 for $\oo_1\subset\oo_\eta$ then yields $\lambda_p(\mathcal J_{A,\oo_1})<\lambda_p(\mathcal J_{A,\oo_\eta})\le s(\mathcal J_{A,\oo_2})$.\vspace{-2mm}\end{proof}

If neither operator admits a principal eigenfunction, the comparison reduces to the maxima of $s(A(x))$ over the two domains.

\begin{coro}\label{cor3.2}
Let $\oo_1\subsetneqq\oo_2$ be bounded connected domains, and
assume that $A$ is irreducible at some point of $\boo_1$. If neither operator admits
a principal eigenvalue, then $s(A,\overline\Omega_2)>s(A,\overline\Omega_1)$ implies $s(\mathcal J_{A,\oo_2})>s(\mathcal J_{A,\oo_1})$, while $s(A,\overline\Omega_2)=s(A,\overline\Omega_1)$ implies $s(\mathcal J_{A,\oo_2})=s(\mathcal J_{A,\oo_1})$. In the second case strict domain monotonicity fails.\vspace{-2mm}
\end{coro}

We next consider intervals $(0,l)$ with a fixed coefficient matrix $A(x)$.\vspace{-1mm}

\begin{coro}\label{cor3.3}
In one dimension, let $l_0>0$, and assume that $A$ is irreducible at some $x_*\in(0,l_0)$. Suppose that $s(A(x))$ is nondecreasing on $[l_0,\infty)$ and
  \[s(A(x))\le s(A(l_0))\;\;\;\text{for all}\;\; 0\le x\le l_0.\]
Then $s(\mathcal J_{A,(0,l)})$ is strictly increasing for $l>l_0$.\vspace{-1mm}
\end{coro} 

\begin{proof} Let $l_2>l_1>l_0$. For $0\le x\le l_0$, we have $s(A(x))\le s(A(l_0))\le s(A(l_i))$; for $l_0\le x\le l_i$, monotonicity implies $s(A(x))\le s(A(l_i))$.
Consequently, $s(A,[0,l_i])=s(A(l_i))$ for $i=1,2$.

If $s(A(l_2))=s(A(l_1))$, monotonicity gives\vspace{-1mm}
  \[s(A(x))=s(A(l_2))=s(A,[0,l_2]),\;\;\;\forall\; x\in[l_1,l_2].\vspace{-1mm}\]
Since the interval $(l_1,l_2)$ has positive length,
Proposition~\ref{prop2.1} implies that $\mathcal J_{A,(0,l_2)}$ admits a principal
eigenvalue. By Theorem~\ref{th1.1}\textup{(ii)}, we obtain
$s\big(\mathcal J_{A,(0,l_1)}\big)<s\big(\mathcal J_{A,(0,l_2)}\big)$.

Now suppose $s(A(l_2))>s(A(l_1))$. If $\mathcal J_{A,(0,l_1)}$ has a principal eigenvalue, apply
Theorem~\ref{th1.1}\textup{(ii)} directly. Otherwise, Theorem~\ref{th2.1}\textup{(iii)} gives $s\big(\mathcal J_{A,(0,l_1)}\big)=s(A,[0,l_1])$, and\vspace{-1mm}
\[s\big(\mathcal J_{A,(0,l_2)}\big)\ge s(A,[0,l_2])=s(A(l_2))>s(A(l_1))
=s(A,[0,l_1])=s\big(\mathcal J_{A,(0,l_1)}\big).\]
The proof is complete.\vspace{-2mm}\end{proof}

\section{Domain limits and critical sizes}\label{sec4}

In this section, we study the small- and large-domain limits  of the spectral bound and use these limits to determine critical domain sizes.

\subsection{Small-domain limits---Proof of Theorem~\ref{th1.2}\textup{(i)}}

\begin{proof}[Proof of Theorem~\ref{th1.2}\textup{(i)}]
Set $\mu_0=s(A(x_0))$. By the Perron--Frobenius theorem, there exists $v\gg0$ such that $A(x_0)v=\mu_0v$. By continuity, there exists $\sigma_0>0$ such that
$a_{ij}(x)>0$ whenever $a_{ij}(x_0)>0$ and $|x-x_0|<\sigma_0$. Thus $A(x)$ is irreducible for $|x-x_0|<\sigma_0$, so the Collatz--Wielandt formulas \eqref{2.2} apply on each
$ B_\sigma(x_0)$ for sufficiently small $\sigma$.

Define
\[
\eta(\sigma)=\max_{i\in\mathcal S}\,\sup_{|x-x_0|\le\sigma}
\frac{|((A(x)-A(x_0))v)_i|}{v_i}\to0,\;\;\;
C=\max_{i\in\mathcal S}d_i\|J_i\|_{L^\infty(\mathbb R^{2n})}.\]
Taking $v$ as a constant test function, we have, for all
$x\in\overline B_\sigma(x_0)$ and $i\in\mathcal S$,
\[\mu_0-\eta(\sigma)\le\frac{\mathcal J_{i, A,B_\sigma(x_0)}[v](x)}{v_i}
\le\mu_0+\eta(\sigma)+C|B_\sigma(x_0)|.\]
By \eqref{2.2},
  \[\mu_0-\eta(\sigma)\le s\big(\mathcal J_{A,B_\sigma(x_0)}\big)
  \le\mu_0+\eta(\sigma)+C|B_\sigma(x_0)|.\]
Both bounds converge to $\mu_0$, which yields the desired limit.
These estimates remain valid on each connected component of $\Omega\cap B_\sigma(x_0)$, with constants uniform across components; hence the convergence is uniform over these components.\www
\end{proof}

\subsection{Large-domain limits---Proof of Theorem~\ref{th1.2}\textup{(ii)}}

%In this proof, write $\Omega_l=B_l(0)$.\vspace{-2mm}

\begin{proof}[Proof of Theorem~\ref{th1.2}\textup{(ii)}]
Since $A^\infty+D$ is a constant cooperative irreducible
matrix, $s(A^\infty+D)$ is the principal eigenvalue of $A^\infty+D$. Choose
$v^\infty\gg0$ with $(A^\infty+D)v^\infty=s(A^\infty+D) v^\infty$.
Note that $\int_{B_l(0)} K_i(y-x)\dy\le1$ and $A(x)\le A^\infty$. The constant test function $v^\infty$ satisfies
 \[\mathcal J_{i, A,B_l(0)}[v^\infty](x)\le d_i v^\infty_i+ A^\infty_{i,*} v^\infty
=s(A^\infty+D) v^\infty_i.\]
The inf-sup formula \eqref{2.2} yields $s\big(\mathcal J_{A,B_l(0)}\big)\le s(A^\infty+D)$ for all sufficiently large $l$ with $x_0\in B_l(0)$.

Choose $\ell_k\to\infty$ sufficiently large that $x_0\in B_{\ell_k}(0)=:B_{\ell_k}$ and $H_k\subset B_{\ell_k}$. For each $i\in\mathcal S$, set
\[
\theta_{i,k}:=\frac1{|H_k|}\int_{H_k}\int_{H_k}K_i(y-x)\dy\dx.
\]

\textit{Claim 1.} For each $i\in\mathcal S$, $\theta_{i,k}\to1$ as $k\to\infty$, and
\[
\frac1{|H_k|}\int_{H_k}
\frac{\int_{B_{\ell_k}}K_i(y-x)\psi_i(y)\dy}{\psi_i(x)}\dx
\ge\theta_{i,k}
\]
for every $\psi\in[C(\overline{B}_{\ell_k})]^m$ with $\psi\gg0$. Indeed, with $z=y-x$, the conditions $x,y\in H_k$ are equivalent to $x\in H_k\cap(H_k-z)$, and $|H_k\cap(H_k-z)|=\prod_{j=1}^n(l_k-|z_j|)_+$. Consequently,
  \[\theta_{i,k}=\int_{\mathbb R^n} K_i(z)\prod_{j=1}^n\left(1-\frac{|z_j|}{l_k}\right)_+\dz
 \to1.\]
The convergence follows from dominated convergence, since the product lies in $[0,1]$, tends pointwise to $1$, and $K_i\in L^1(\mathbb R^n)$.

Restrict the integral to $H_k$ and use the symmetry of $K_i$.
The inequality $t+t^{-1}\ge2$ with $t=\psi_i(y)/\psi_i(x)$ then gives
 \begin{equation}\label{4.3}
\begin{aligned}[b]
\frac1{|H_k|}\int_{H_k}\frac{\int_{B_{\ell_k}} K_i(y-x)\psi_i(y)\dy}{\psi_i(x)}\dx
&\ge\frac1{|H_k|}\int_{H_k}\int_{H_k} K_i(y-x)\frac{\psi_i(y)}{\psi_i(x)}\dy\dx\\
&=\frac1{2|H_k|}\int_{H_k}\int_{H_k} K_i(y-x)
 \left(\frac{\psi_i(y)}{\psi_i(x)}+\frac{\psi_i(x)}{\psi_i(y)}\right)\dy\dx\\
&\ge\frac1{|H_k|}\int_{H_k}\int_{H_k} K_i(y-x)\dy\dx
 =\theta_{i,k}.
\end{aligned}\zzz
\end{equation}
Claim 1 is proved. 
    
Set\zzz
\[
\delta_k=\max_{1\le i,j\le m}\,\sup_{x\in H_k}|a_{ij}(x)-a^\infty_{ij}|,
\;\;\; b_{ii,k}=a^\infty_{ii}-\delta_k,\;\;\;
b_{ij,k}=(a^\infty_{ij}-\delta_k)_+\;\;\;(i\ne j),\www
\]
and write $B_k=(b_{ij,k})$. Then $B_k\le A(x)$ on $H_k$,
$B_k\to A^\infty$, and $B_k$ is cooperative. Every positive
off-diagonal entry of $A^\infty$ remains positive in $B_k$ for all
sufficiently large $k$. Hence $B_k$ is irreducible for such indices.

\textit{Claim 2.} For every sufficiently large $k$, there exist weights
$p_{i,k}>0$, with $\sum_i p_{i,k}=1$, such that
  \[\sum_i p_{i,k}\frac{(B_k u)_i}{u_i}\ge s(B_k+D)-\sum_i p_{i,k}d_i\]
for every vector $u\gg0$. Indeed, take $v_k,w_k\gg0$ to be a right and a left eigenvector, respectively, satisfying  \zzz
  \[ (B_k+D)v_k=s(B_k+D)v_k,\;\;\; w_k^{\mathsf T}(B_k+D)=s(B_k+D)w_k^{\mathsf T},\zzz\]
normalized by $w_k^{\mathsf T}v_k=1$. Set $p_{i,k}=w_{i,k}v_{i,k}>0$ and
$q_{ij,k}=w_{i,k}b_{ij,k}v_{j,k}\ge0$ for $i\ne j$.
The normalization implies $\sum_i p_{i,k}=1$.
Multiplying row $i$ of the right eigenvector equation by $w_{i,k}$
and column $i$ of the left eigenvector equation by $v_{i,k}$, we obtain
\begin{equation}\label{4.4}
\sum_{j\ne i}q_{ij,k}=(s(B_k+D)-b_{ii,k}-d_i)p_{i,k} =\sum_{j\ne i}q_{ji,k}.
\end{equation}
For any vector $u\gg0$, set $g_i=\log(u_i/v_{i,k})$. Then\zzz
\[
\sum_{i\ne j}q_{ij,k}(g_j-g_i)
=\sum_i g_i\bigg(\sum_{j\ne i}q_{ji,k}-\sum_{j\ne i}q_{ij,k}\bigg)=0.\zzz
\]
Using \eqref{4.4} and ${\rm e}^{g_j-g_i}\ge1+g_j-g_i$, we obtain\zzz
 \begin{align}
\sum_i p_{i,k}\frac{(B_k u)_i}{u_i}&=\sum_i p_{i,k}b_{ii,k}
  +\sum_{i\ne j}q_{ij,k} {\rm e}^{g_j-g_i}\nonumber\\
&\ge\sum_i p_{i,k}b_{ii,k}+\sum_{i\ne j}q_{ij,k}
 =s(B_k+D)-\sum_i p_{i,k}d_i.\label{4.5}\www
\end{align}
Claim 2 is proved.

Integrating \eqref{4.5} with $u=\psi(x)$ over $H_k$ and dividing by
$|H_k|$, we add \eqref{4.3} multiplied by $p_{i,k}d_i$ and sum over $i$. Since $B_k\le A(x)$ on $H_k$, $p_{i,k}>0$, and $\sum_i p_{i,k}=1$, for every $\psi\in[C(\overline{B_{\ell_k}})]^m$ with $\psi\gg0$ we have
\begin{align*}
\sup_{\substack{x\in B_{\ell_k},\,i\in\mathcal S}}
 \frac{\mathcal J_{i, A,B_{\ell_k}}[\psi](x)}{\psi_i(x)}
&\ge\frac1{|H_k|}\sum_i p_{i,k}\int_{H_k}
 \frac{\mathcal J_{i, A,B_{\ell_k}}[\psi](x)}{\psi_i(x)}\dx\\
&\ge s(B_k+D)-\sum_i p_{i,k} d_i(1-\theta_{i,k}).
\end{align*}
Since $B_k\to A^\infty$, we have $s(B_k+D)\to s(A^\infty+D)$.
Furthermore,
\[
0\le\sum_i p_{i,k}d_i(1-\theta_{i,k})
\le\max_i d_i\,\max_i(1-\theta_{i,k})\to0.\vspace{-2mm}
\]
Taking the infimum over $\psi$ and using \eqref{2.2}, we obtain
\begin{equation}\label{4.6}
  s(\mathcal J_{A,B_{\ell_k}})\ge s(B_k+D)-\sum_i p_{i,k}d_i(1-\theta_{i,k})\to s(A^\infty+D)\;\;\;\text{as}\;\;k\to\infty.\vspace{-2mm}
\end{equation}
For all sufficiently large $l$, the ball $B_l(0)$ contains $x_0$. By \eqref{3.2}, $ s(\mathcal J_{A,B_l(0)})$ is nondecreasing for such $l$. Combining this monotonicity with the upper bound and \eqref{4.6}, we obtain \eqref{1.4}.\www
\end{proof}

The next two examples show that \eqref{1.2} cannot be omitted and
that pointwise convergence along a sequence cannot replace \eqref{1.3}.\zzz

\begin{example}\label{ex4.1}
The bound $A(x)\le A^\infty$ in Theorem~\ref{th1.2}\textup{(ii)} cannot be omitted.
 For example, take $n=m=1$, $J(x,y)=(1-|x-y|)_+$, $d>0$, $L>0$, and $a_0>d$. Define
\[
a(x)=a_0\left(1-{|x|}/{L}\right)_+,
\;\;\; a^\infty=0.
\]
All the other hypotheses of the theorem are satisfied, but the asserted
limit does not hold.
\end{example}

\begin{proof}
Since $a=0$ outside $[-L,L]$, condition \eqref{1.3} holds on
$H_k=[L+k,L+2k]$. However, Theorem~\ref{th2.1} implies
 \[s\big(\mathcal J_{a,B_l(0)}\big)\ge\max_{x\in\overline{B}_l(0)}a(x)
 =a_0>d=s(a^\infty+d),\;\;\;\forall\;l>0.\]
Thus \eqref{1.4} cannot hold.\vspace{-2mm}
\end{proof}

\begin{example}\label{ex4.2}
Convergence $A(x_k)\to A^\infty$ along a sequence $x_k\to\infty$
cannot replace \eqref{1.3} in Theorem~\ref{th1.2}\textup{(ii)}.
 For example, take $n=m=1$, $J(x,y)=K(y-x)$ and
\[
d=1,\;\;\; K(z)=(1-|z|)_+,\;\;\;
a(x)=2+\cos(2\pi x),\;\;\; a^\infty=3.\]
Then $1\le a(x)\le a^\infty$ and $a(k)=a^\infty$ for every positive
integer $k$, but
\begin{equation}\label{4.7}
3\le\lim_{l\to\infty} s\big(\mathcal J_{a,B_l(0)}\big)\le 11/3<4=a^\infty+1.\zzz
\end{equation}
\end{example}

\begin{proof} Take $\varphi=a$. Since $K$ is even and $\int_{-1}^{1}(1-|z|)\cos(2\pi z)\dz=0$, we have $\int_{\mathbb R}K(y-x)\varphi(y)\dy=2$. Since $1\le\varphi\le3$, its restriction to each 
$\overline{B}_l(0)$ is a strongly positive test function, and
  \[\frac{\mathcal J_{a,B_l(0)}[\varphi](x)}{\varphi(x)}
\le a(x)+\frac2{a(x)}\le\max_{1\le t\le3}\left(t+\frac2t\right)=\frac{11}{3}.\]
The inf-sup formula \eqref{2.2} and the lower bound \eqref{2.4} yield
\[3=s(a,\overline{B}_l(0))\le s\big(\mathcal J_{a,B_l(0)}\big)\le 11/3,\;\;l>0.\]
By \eqref{3.2}, $s\big(\mathcal J_{a,B_l(0)}\big)$ is nondecreasing in $l$, and the conclusion \eqref{4.7} follows.
\end{proof}

\subsection{Critical sizes}

\begin{proof}[Proof of Corollary~\ref{cor1.1}]
 Since $A(0)$ is irreducible and $\mathcal J_{A,B_L(0)}$ has principal eigenvalue zero,
by Theorem~\ref{th1.1}\textup{(ii)}, we have\zzz
\[
s\big(\mathcal J_{A,B_l(0)}\big)<s\big(\mathcal J_{A,B_L(0)}\big)=0\;\;\;\text{for}\;\;0<l<L,\;\;\;
s\big(\mathcal J_{A,B_l(0)}\big)>s\big(\mathcal J_{A,B_L(0)}\big)=0\;\;\;\text{for}\;\;l>L.\zzz
\]
This proves \eqref{1.5} and the uniqueness of the critical size.
Under the hypotheses of Theorem~\ref{th1.2}\textup{(ii)}, the limit \eqref{1.4} and domain monotonicity yield\zzz
  \[s(A^\infty+D)=\lim_{l\to\infty}s\big(\mathcal J_{A,B_l(0)}\big)\ge s\big(\mathcal J_{A,B_{2L}(0)}\big)>0.\zzz\]
This proves the last assertion.\vspace{-2mm}
\end{proof}

The following example has symmetric kernels and coefficients, but its spectral bound vanishes for every $l\ge1$.\zzz

\begin{example}\label{ex4.3} Take $n=1$, $m=2$, and \zzz
 \begin{align*} 
 d_1&=d_2=1,\;\;\; J_1(x,y)=J_2(x,y)=K(y-x):=(1-|x-y|)_+,\\
 b(x)&=4\min\{\sqrt{|x-1|},1\},\;\;\;
A(x)=\begin{pmatrix}-1-b(x)&1\\ 
 1&-1-b(x)\end{pmatrix}.\zzz\end{align*}
The kernels satisfy \textbf{\bf(J)}, and $A$ satisfies \textbf{\bf(A)} and is
irreducible at every point. We have\zzz
\begin{equation}\label{4.8}
 s\big(\mathcal J_{A,B_l(0)}\big)<0\;\;\;\text{for}\;\;0<l< 1,\;\;\;  s\big(\mathcal J_{A,B_l(0)}\big)=0\;\;\;\text{for}\;\;l\ge 1.\zzz
\end{equation}
For every $l\ge1$, the operator has no principal eigenfunction. The conclusions in \eqref{4.8} also hold with $A$ replaced by\zzz
  \[A_l(x)=A(x)-\frac{b(x)}{4l}I.\zzz\]
This family satisfies $A_{l_1}\le A_{l_2}$ and $A_{l_1}\not\equiv A_{l_2}$ whenever $l_1<l_2$.\zzz
\end{example}

\begin{proof} Set $\mathbf e=(1,1)^{\mathsf T}$ and $\beta_l=\min_{\overline{B}_l(0)}b$. We first show that for every $l>0$,
    \begin{equation}\label{4.9}
 -\beta_l=s(A,\overline{B}_l(0))\le s\big(\mathcal J_{A,B_l(0)}\big)\le-\beta_l/4.
    \end{equation}
In fact, since $A(x)\mathbf e=-b(x)\mathbf e$, we have $-\beta_l=s(A,\overline{B}_l(0))$. The left-hand inequality in \eqref{4.9} follows from \eqref{2.4}. Using $0\le K\le1$ and $\int_{\mathbb R}K=1$, we obtain
 \[\int_{B_l(0)}\frac{K(y-x)}{b(y)}\dy\le \frac14\int_{\mathbb R}K(y-x)\dy +\frac14\int_{|y-1|<1}
 \bigl(|y-1|^{-1/2}-1\bigr)\dy=\frac14+\frac12=\frac34.\]
For $\ep>0$, take $\Phi_\ep=\mathbf e/(\ep+b)\gg0$. For $i=1,2$, direct substitution gives
 \[\frac{\mathcal J_{i, A,B_l(0)}[\Phi_\ep](x)}{\Phi_{\ep,i}(x)}
=(\ep+b(x))\int_{-l}^l\frac{K(y-x)}{\ep+b(y)}\dy-b(x)
 \le\frac34(\ep+b(x))-b(x)=\frac34\ep-\frac14b(x). \]
Letting $\varepsilon\searrow0$ and using the inf-sup characterization \eqref{2.2}, we obtain the right-hand inequality in \eqref{4.9}.
        
Since $\beta_l=4\sqrt{1-l}>0$ for $0<l<1$ and $\beta_l=0$ for $l\ge1$, \eqref{4.8} follows. For $l\ge1$, \eqref{4.9} and Theorem~\ref{th2.1} exclude a strongly positive principal eigenfunction. 
        
Finally, since $A_l\le A$, coefficient monotonicity and \eqref{2.4} yield
 \[-[1+1/(4l)]\beta_l=s(A_l,\overline{B_l(0)})
 \le s\big(\mathcal J_{A_l,B_l(0)}\big)\le s\big(\mathcal J_{A,B_l(0)}\big).\]
Hence \eqref{4.8} also holds with $A$ replaced by $A_l$.\www
\end{proof}

The following example shows that evenness in Theorem~\ref{th1.2}\textup{(ii)} cannot simply be omitted.

\begin{example}\label{ex4.4}
Let $n=1$, $\sigma>0$, $b\in\mathbb R$, and
$J_i(x,y)=K_b(y-x)$ for every $i\in\mathcal S$, where
\[
k_\sigma(z)=\frac{{\rm e}^{-z^2/(2\sigma^2)}}{\sqrt{2\pi}\sigma},\;\;\;
K_b(z)=k_\sigma(z-b),\;\;\; c_b={\rm e}^{-b^2/(2\sigma^2)}.
\]
The kernels satisfy \textbf{\bf(J)} and are not even when $b\ne0$.
Write $s_b\big(\mathcal J_{A,B_l(0)}\big)$ for the spectral bound of this family. Assume \textbf{\bf(A)}, and let $A^\infty$ be a constant cooperative irreducible matrix satisfying
\eqref{1.2}--\eqref{1.3}. Then
\begin{equation}\label{4.10}
\lim_{l\to\infty}s_b\big(\mathcal J_{A,B_l(0)}\big)=s(A^\infty+c_bD).
\end{equation}
For fixed $\sigma$, this limit decreases strictly with $|b|$, equals
$s(A^\infty+D)$ at $b=0$, and tends to $s(A^\infty)$ as $|b|\to\infty$.
In the scalar case with constant coefficient $a$, the limit is
$a+dc_b$. In particular, for $a=0$ and $d=1$, the limit is $c_b<1$ for $b\ne0$,
although $\int_{\mathbb R}K_b=1$. Thus evenness in
Theorem~\ref{th1.2}\textup{(ii)} cannot simply be omitted.\www
\end{example}

\begin{proof} Define $S_b[\varphi](x)={\rm e}^{-bx/\sigma^2}\varphi(x)$. Then both $S_b$ and $S_b^{-1}$ are bounded on $[C(\overline{B}_l(0))]^m$. The identity ${\rm e}^{-bz/\sigma^2}K_b(z)=c_bk_\sigma(z)$ yields
 \[ \bigl(S_b^{-1}\mathcal J_{A,B_l(0)}S_b[\varphi]\bigr)_i(x)
  =c_bd_i\int_{B_l(0)}k_\sigma(y-x)\varphi_i(y)\dy +A_{i,*}(x)\varphi(x),\]
so the two operators are similar and have the same spectrum.

The transformed operator has an even normalized kernel, integral coefficients $c_bd_i$, and the unchanged local matrix $A(x)$. Theorem~\ref{th1.2}\textup{(ii)} yields \eqref{4.10}. The assertions about $b$ follow from $c_b={\rm e}^{-b^2/(2\sigma^2)}$ and the continuity and strict monotonicity of the spectral bound of cooperative irreducible matrices, since every $d_i>0$.
    
The exponent $-b/\sigma^2$ is the unique minimizer of the moment function:
  \[\int_{\mathbb R}K_b(z){\rm e}^{pz}\dz={\rm e}^{bp+\sigma^2p^2/2}
  =c_b\exp\!\left[\frac{\sigma^2}{2}\left(p+\frac b{\sigma^2}\right)^2\right],\]
and the corresponding minimum value of this function is $c_b$. This similarity transformation is used only on bounded intervals: for $b\ne0$, neither exponential multiplier is bounded on $\mathbb R$.
\end{proof}

\section{An application to a heterogeneous epidemic model}\label{sec5}

In the rest of the paper, for intervals $(0,l)$, we use the notation
    \[\mathcal J_{a,l}:=\mathcal J_{a,(0,l)},\;\;\; s(\mathcal J_{a,l})=s\big(\mathcal J_{a,(0,l)}\big).\]
The same convention applies to matrix coefficients. Throughout this section, we fix $n=1$, $m=2$, and $D=\operatorname{diag}(d_1,d_2)$.
    
The preceding spectral results allow us to determine how habitat size dictates the extinction or persistence of the epidemic. Motivated by Li, Li, and Wang~\cite{LLW25SAM}, we consider the following heterogeneous system posed on the fixed interval $(0,l)$:
\begin{equation}\label{5.1}
 \begin{cases}
 u_t=d_1\dd\int_0^l  K_1(y-x)u(t,y)\dy-(d_1+a(x))u+H(x,v),\\[2mm]
 v_t=d_2\dd\int_0^l  K_2(y-x)v(t,y)\dy -(d_2+b(x))v+G(x,u),
\end{cases}
\end{equation}
for $t>0$ and $x\in[0,l]$, with initial data
$(u_0,v_0)\ge0$ in $[C([0,l])]^2$.
The dispersal terms correspond to zero density outside $(0,l)$.

Assume \textbf{(J)} holds with $J_i(x,y)= K_i(y-x)$ and $ K_i$ even,
for $i=1,2$. Suppose that $H(x,0)=G(x,0)=0$ and that $H,G$ are
differentiable in their second variable at zero. Set
$h(x)=H_v(x,0)$ and $g(x)=G_u(x,0)$. Then the linearized operator
at zero is $\mathcal J_{A,l}$, where
 \[A(x)=\begin{pmatrix}-d_1-a(x)&h(x)\\ 
 g(x)&-d_2-b(x)\end{pmatrix}.\]
We impose the following assumptions:\zzz
\begin{enumerate}[leftmargin=14mm]
\item[\textbf{(E1)}] The functions $a,b,h,g$ are continuous, bounded, and positive
on $\mathbb R_+$.\zzz
\item[\textbf{(E2)}] The coefficients satisfy $(d_1+a(0))(d_2+b(0))>h(0)g(0)$, and
\begin{equation}\label{5.2}
\text{there exists }L_0>0\text{ such that }s(\mathcal J_{A,L_0})>0.\www
\end{equation}
\end{enumerate}

Next, we show the threshold dynamics of \eqref{5.1} in terms of the interval length $l$.\www

\begin{prop}\label{prop5.1} Assume \textbf{\bf(E1)} and \textbf{\bf(E2)} hold.
Then $s(\mathcal J_{A,l})$ is continuous and nondecreasing in $l>0$, and
  \begin{equation}\label{5.3}
\lim_{l\searrow0}s(\mathcal J_{A,l})=s(A(0))<0.
  \end{equation}
The number
  \[l_*:=\inf\{l>0:s(\mathcal J_{A,l})>0\}=\sup\{l>0:s(\mathcal J_{A,l})\leq0\}\]
belongs to $(0,\infty)$ and satisfies 
\begin{equation}\label{5.4}
 s(\mathcal J_{A,l_*})=0,\;\;\; s(\mathcal J_{A,l})\le0\;\;\;\text{for}\;\; 0<l\le l_*,\;\;\;
s(\mathcal J_{A,l})>0\;\;\;\text{for}\;\;l>l_*.
\end{equation}
The set $\{l>0:s(\mathcal J_{A,l})=0\}$ is a closed interval with
right endpoint $l_*$, possibly a singleton.\zzz
\end{prop}

\begin{proof} By \textbf{(E1)}, $A(x)$ is irreducible for every $x\ge0$.
Proposition~\ref{prop3.1} and \eqref{3.2} imply that
$l\mapsto s(\mathcal J_{A,l})$ is continuous and nondecreasing.
For decreasing intervals, take the maps $x\mapsto\min\{x,l\}$ in
Proposition~\ref{prop3.1}\textup{(ii)}. As in the proof of Theorem~\ref{th1.2}\textup{(i)}, the limit in \eqref{5.3} holds. Since $A(0)$ has negative trace and
$\det A(0)=(d_1+a(0))(d_2+b(0))-h(0)g(0)>0$, we have $s(A(0))<0$.

By \eqref{5.2} and \eqref{5.3}, $0<l_*<\infty$. The definition of $l_*$ and monotonicity imply that
$s(\mathcal J_{A,l})\le0$ for $l<l_*$ and $s(\mathcal J_{A,l})>0$ for $l>l_*$. Continuity yields
$s(\mathcal J_{A,l_*})=0$ and hence \eqref{5.4}. The last assertion follows from continuity and monotonicity.\vspace{-2mm}
\end{proof}

Condition \eqref{5.2} holds if $(a(x),b(x),h(x),g(x))\to(a^\infty,b^\infty,h^\infty,g^\infty)$ as $x\to\infty$, where $a^\infty,b^\infty\ge0$, $h^\infty,g^\infty>0$, and $h^\infty g^\infty>a^\infty b^\infty$. Indeed, set $A^\infty=\lim_{x\to\infty}A(x)$ and
 \[ B_\eta=\begin{pmatrix}
-d_1-a^\infty-\eta&h^\infty-\eta\\ 
g^\infty-\eta&-d_2-b^\infty-\eta
 \end{pmatrix}.\]
For all sufficiently small $\eta>0$, we have $s(B_\eta+D)>0$ and $B_\eta\le A(x)$ for $x\ge R_\eta$,  with some $R_\eta>0$. By Lemma~\ref{lem3.2}\textup{(i)}, \eqref{3.2}, and Theorem~\ref{th1.2}\textup{(ii)}, for $l>R_\eta$,
 \[s(\mathcal J_{A,l})\ge s\big(\mathcal J_{B_\eta,(R_\eta,l)}\big)
=s\big(\mathcal J_{B_\eta,l-R_\eta}\big)\to s(B_\eta+D)>0\;\;\;\text{as }l\to\infty,\]
where the equality follows by translation. Since $s(\mathcal J_{A,l})$ is nondecreasing and bounded above by $\sup_{x\ge0}\max\{h(x)-a(x),g(x)-b(x)\}$, by \eqref{2.2}, letting $\eta\searrow0$ gives
 \[\lim_{l\to\infty}s(\mathcal J_{A,l})\ge s(A^\infty+D)>0,\]
which proves \eqref{5.2}. If $A(x)\le A^\infty$ for all $x\ge0$, testing \eqref{2.2} with a positive Perron vector of $A^\infty+D$ gives the reverse inequality, so equality holds.
    
More generally, the same comparison proves \eqref{5.2} if the four coefficients converge uniformly to these constants on intervals $[r_k,r_k+l_k]$, where $r_k>0$ and $l_k\to\infty$. Indeed, for fixed small $\eta>0$ and all sufficiently large $k$,
 \[ s(\mathcal J_{A,r_k+l_k}) \ge s\big(\mathcal J_{B_\eta,(r_k,r_k+l_k)}\big)
 =s\big(\mathcal J_{B_\eta,l_k}\big)\to s(B_\eta+D)>0.\]
    
If, in addition to \textbf{(E1)} and \textbf{(E2)}, $a,b$ are
nonincreasing and $h,g$ are nondecreasing, then $s(A(x))$ is
nondecreasing. Corollary~\ref{cor3.3} implies that
$l\mapsto s(\mathcal J_{A,l})$ is strictly increasing.
In this case, $s(\mathcal J_{A,l})<0$ for $0<l<l_*$, and $l_*$
is the unique zero of the spectral bound.

To apply \cite[Theorem B]{WZ26JMPA} to the problem \eqref{5.1}, we impose the following additional assumptions:\zzz
\begin{enumerate}[leftmargin=14mm]
\item[\textbf{(N1)}] The functions $H,G$ and their derivatives in the
second variable are continuous on $\mathbb R_+\times\mathbb R_+$, and
$H_v(x,v)>0, G_u(x,u)>0$. \zzz
\item[\textbf{(N2)}] For each $x\ge0$, both $H(x,z)/z$ and
$G(x,z)/z$ are strictly decreasing for $z>0$.\zzz
\item[\textbf{(N3)}] For each $l>0$, there exist constants $P_l,Q_l>0$
such that
  \[H(x,Q_l)\le a(x)P_l,\;\;\; G(x,P_l)\le b(x)Q_l,\;\;\; x\in[0,l].\zzz\]
\end{enumerate}
Condition \textbf{(N2)} follows if $H$ and $G$ are strictly concave
in their second variable, as assumed in~\cite{LLW25SAM}.
Condition \textbf{(N3)} holds if $H(x,z)/z$ and $G(x,z)/z$ tend to
zero as $z\to\infty$, uniformly on each bounded interval in $x$.

Under \textbf{(E1)} and \textbf{(N1)--(N3)}, for any given $l>0$ and any nonnegative initial datum in $[C([0,l])]^2$, problem \eqref{5.1} admits a unique bounded global solution $(u,v)$. 
The sign of the spectral bound determines the limiting behavior as follows.\zzz

\begin{coro}\label{cor5.1}
Assume \textbf{\bf(E1)}, \textbf{\bf(E2)} and \textbf{\bf(N1)--(N3)} hold.
Let $l_*$ be given by Proposition~\ref{prop5.1}. \www
\begin{enumerate}[label=\textup{(\roman*)},leftmargin=8mm]
\item If $0<l\le l_*$, zero is the unique nonnegative equilibrium
and is globally asymptotically stable:
 \[\lim_{t\to\infty}\bigl(\|u(t,\cdot)\|_\infty+\| v(t,\cdot)\|_\infty\bigr)=0.\zzz\]
Moreover, if $s(\mathcal J_{A,l})<0$, the convergence is exponential.\www
\item If $l>l_*$, then the zero equilibrium is unstable, and \eqref{5.1} has a unique positive equilibrium $(U_l,V_l)\gg0$ that is globally asymptotically stable for
nontrivial nonnegative initial data. In particular, for every $(u_0,v_0)\not\equiv(0,0)$,
 \[\lim_{t\to\infty}\bigl(\|u(t,\cdot)-U_l\|_\infty+\|v(t,\cdot)-V_l\|_\infty\bigr)=0.\www\]
\end{enumerate}
\end{coro}

\begin{proof}    Set
 \[f(x,u,v)=(-a(x)u+H(x,v),\;\;-b(x)v+G(x,u))^{\mathsf T}.\]
By \textbf{(N1)}, $f$ is continuously differentiable in $(u,v)$, cooperative, and its Jacobian is irreducible for $u,v\ge0$. By \textbf{(N2)}, for $u,v>0$ and $0<\rho<1$, we have
 \[ f(x,\rho u,\rho v)-\rho f(x,u,v)=\begin{pmatrix}
  H(x,\rho v)-\rho H(x,v)\\
   G(x,\rho u)-\rho G(x,u)
   \end{pmatrix}\gg0.\]
Thus $f$ is strongly subhomogeneous. Since $\int_0^l K_i(y-x)\dy\le1$, \textbf{(N3)} implies that $(P_l,Q_l)$ is a positive upper solution. By \textbf{(N2)}, $c(P_l,Q_l)$ is also an upper solution for every $c\ge1$. With $c$ sufficiently large, the comparison principle ensures a bounded global solution for each nonnegative initial datum. By \textbf{(J)} and \textbf{(N1)}, every nontrivial nonnegative solution is strongly positive for $t>0$.
    
By \cite[Theorem B(1)--(3)]{WZ26JMPA} and \eqref{5.4}, we obtain the stated conclusions on equilibria and convergence. Since every nontrivial nonnegative solution is strongly positive for $t>0$, the convergence holds for all nontrivial nonnegative initial data, and every nontrivial nonnegative equilibrium is strongly positive.
    
For $0<l\le l_*$, stability of the zero equilibrium follows from continuous dependence and comparison with the solution starting from a fixed upper solution $\overline w\gg0$, which converges to zero. For $l>l_*$, set $w=(u,v)^{\mathsf T}$ and $w_l=(U_l,V_l)^{\mathsf T}\gg0$. For $0<\rho<1$, \textbf{(N2)} implies that $\rho w_l$ and $\rho^{-1}w_l$ are a lower and an upper solution, respectively. By the comparison principle,
  \[\rho w_l\le w(0,\cdot)\le\rho^{-1}w_l\implies
\rho w_l\le w(t,\cdot)\le\rho^{-1}w_l\;\;\;(t\ge0).\]
The two bounds converge uniformly to $w_l$ as $\rho\nearrow1$, which implies stability of $w_l$. The trivial equilibrium is unstable, as every nontrivial nonnegative solution converges to $w_l$.\www
\end{proof}

\begin{remark}[An interval of lengths with zero spectral bound]\label{rem5.1}
 Corollary~\ref{cor5.1} shows that a unique extinction--persistence
threshold exists even when the spectral bound vanishes on an interval
of lengths. Since extinction also occurs throughout this interval,
the threshold $l_*$ is its right endpoint.
No principal eigenfunction at zero is required.\www
\end{remark}

\section{A scalar free-boundary problem}
\label{sec6}

We consider a logistic equation with a free boundary and a sign-changing growth rate:
\begin{equation}\label{6.3}
\begin{cases}
\displaystyle
u_t = d\int_0^{h(t)}\!K(y-x)u(t,y)\dy + \big(m(x)-d\big)u - u^2,
& t>0,\ 0\le x<h(t),\\[2mm]
\displaystyle
h'(t) = \mu\int_0^{h(t)}\int_{h(t)}^\infty K(y-x)u(t,x)\dy\dx,& t>0,\\[.2mm]
h(0)=h_0,\;\;\; u(t,h(t))=0,& t>0,\\
u(0,x)=u_0(x),& x\in[0,h_0].
\end{cases}
\end{equation}
Here the kernel $K$ satisfies
\begin{equation}\label{6.1}
K\in C(\mathbb R)\cap L^\infty(\mathbb R),\;\;\; K\ge 0,\;\;\; K(0)>0,\;\;\;
K(z)=K(-z),\;\;\; \int_{\mathbb R}K(z)\,\dz=1,
\end{equation}
and the growth function $m(x)$ is allowed to change sign. We further assume
\begin{equation}\label{6.2}
m\in C_b(\mathbb R_+),\;\;\; d,\mu,h_0>0,\;\;\;
u_0\in C([0,h_0]),\;\;\; u_0>0 \text{ on } [0,h_0),\;\;\; u_0(h_0)=0.
\end{equation}
The density $u$ is understood to vanish outside $[0,h(t)]$.
Since $K$ is even, the outgoing kernel $K (x-y)$ equals $K (y-x)$,
as used in the equation for $h'$.
The corresponding local diffusion problem was investigated by Wang~\cite{W15JDE}.

Set $a(x)=m(x)-d$. For a bounded interval $I\subset(0,\infty)$, write
\[
\mathcal K_I[\varphi](x)=\int_I K(y-x)\varphi(y)\dy,
\;\;\;
\mathcal J_{a,I}[\varphi]=d\mathcal K_I[\varphi]+a\varphi.
\]
Thus $s(\mathcal J_{a,I})$ has the same meaning as in Section~\ref{sec2}; the coefficient 
$a$ includes the loss term $-d$. 

The reaction term $ u(m(x)-u)$ is continuous, locally Lipschitz in $u$ uniformly in $x$, and negative for large $u$, so it satisfies the hypotheses of the well-posedness and comparison results in \cite{LLW22JDE}. Consequently, \eqref{6.3} has a unique global solution, which is increasing with respect to $\mu$, and
\begin{equation*}
0<u(t,x)\le M_0:=\max\kk\{\|u_0\|_\infty,\;\|m\|_\infty\rr\},
\;\;\; 0<h'(t)\le\mu M_0h(t),
\end{equation*}
for $t>0$ and $0\le x<h(t)$. The same upper bound holds at $t=0$.
Continuous dependence on $\mu$ on bounded time intervals will be proved in Lemma~\ref{lem6.2}.

By \cite[Theorem B]{WZ26JMPA},  if $s(\mathcal J_{a,l})>0$, then the stationary problem 
\bee\mathcal J_{a,l}[U]-U^2=0,\;\;\;0\le x\le l\label{6.4}\eee 
has a unique continuous positive solution $U_l$, which is globally asymptotically stable.
If $s(\mathcal J_{a,l})\le0$, every nonnegative solution converges uniformly to zero;
the convergence is exponential when $s(\mathcal J_{a,l})<0$. 

We say that spreading occurs if $h_\infty:=\lim_{t\to\infty}h(t)=\infty$,
and that vanishing occurs if $h_\infty<\infty$ and
$\|u(t,\cdot)\|_{C([0,h(t)])}\to0$.

\subsection{A spreading--vanishing threshold}

Assume that the spectral bound is positive on at least one interval:
\begin{equation}\label{6.5}
\text{there exists }l_0>0\text{ such that }s(\mathcal J_{a, l_0})>0.
\end{equation}
By Proposition~\ref{prop3.1} and \eqref{3.2}, $s(\mathcal J_{a,l})$ is continuous and nondecreasing in $l$, and $s_\infty:=\lim_{l\to\infty}s(\mathcal J_{a,l})$ exists and is finite since  $s(\mathcal J_{a,l})\le\sup m$. No convergence of $m(x)$ at infinity is required. Define
 \bee 
  l_*:=\inf\{l>0:s(\mathcal J_{a,l})>0\}.\label{6.5a}\eee
Condition \eqref{6.5} ensures $l_*<\infty$.
We allow $l_*=0$. If $l_*>0$, continuity implies
 \[s(\mathcal J_{a,l})\le0\;\;\;\text{if}\;\;0<l\le l_*,\;\;\;
s(\mathcal J_{a,l_*})=0,\;\;\; s(\mathcal J_{a,l})>0\;\;\;\text{if}\;\;l>l_*.\zzz\]

\begin{theo}[A spreading--vanishing threshold]\label{th6.1}
Under \eqref{6.1}, \eqref{6.2}, and \eqref{6.5}, let $(u_\mu,h_\mu)$ solve \eqref{6.3} and write $h_\mu(\infty)=\lim_{t\to\infty}h_\mu(t)$.\zzz
\begin{enumerate}[label=\textup{(\roman*)}]
\item If $h_\mu(\infty)<\infty$, then
 \[s(\mathcal J_{a,h_\mu(\infty)})\le0,\;\; h_\mu(\infty)\le l_*,\;\;
 \lim_{t\to\infty}\|u_\mu(t,\cdot)\|_{C([0,h_\mu(t)])}=0.\zzz\]
\item If $h_0\ge l_*$, then $h_\mu(\infty)=\infty$ for every $\mu>0$.\zzz
\item If $0<h_0<l_*$, there exists $\mu^*\in(0,\infty)$ such that
  \[\begin{cases}
h_\mu(\infty)<\infty,\;\;\lim_{t\to\infty}\|u_\mu(t,\cdot)\|_\infty=0,
   &0<\mu\le\mu^*,\\
  h_\mu(\infty)=\infty,&\mu>\mu^*.\end{cases}\]
At the critical parameter,
\begin{equation}\label{6.6}
 h_{\mu^*}(\infty)=l_*.\zzz
\end{equation}
\end{enumerate}
We do not assume that zero is the principal eigenvalue or that
$s(\mathcal J_{a,l})$ is strictly increasing in $l$.\www
\end{theo}

\begin{remark}\label{rem6.1}
As in Remark~\ref{rem5.1}, $l_*$ is the right endpoint of the zero
spectral interval whenever this interval is nonempty.
If $s(\mathcal J_{a,h_0})=0$ and $h_0<l_*$, Theorem~\ref{th6.1}
still gives a threshold $\mu^*\in(0,\infty)$: vanishing occurs for
$0<\mu\le\mu^*$ and spreading for $\mu>\mu^*$.
If $h_0=l_*$, spreading occurs for every $\mu>0$.\zzz
\end{remark}

\begin{lem}\label{lem6.1}
Let $0<l<\ell$ and $s(\mathcal J_{a,\ell})\le0$. Then there exists a constant $C_l$ such that for every $g\in C([0,l])$ with $g\ge0$, the solution of 
  \[v_t=\mathcal J_{a,l}[v],\;\;x\in[0,l],\;t>0;\;\;\;  v(0,x)=g(x),\;\;x\in[0,l]\]
satisfies
\begin{equation}\label{6.7}
\int_0^\infty\!\!\int_0^l v(t,x)\dx\dt\le C_l\|g\|_\infty<\infty.
\end{equation}
 \end{lem}

\begin{proof}\textit{Step 1.}
Put $b(x)=-a(x)$. Since $\max_{[0,\ell]}a(x)\le s(\mathcal J_{a,\ell})\le0$, one has $b\ge0$ on $[0,\ell]$.
We first prove
\begin{equation}\label{6.8}
1/b\in L^1(0,\ell).
\end{equation}
For each $\ep>0$, by the inf-sup formula in \eqref{2.2}, there exists $\varphi_\ep\gg0$ such that
  \bee
d\int_0^\ell K(y-x)\varphi_\ep(y)\dy\le (-a(x)+\ep)\varphi_\ep(x)=(b(x)+\ep)\varphi_\ep(x).\label{6.9}\eee
Choose $\delta,\kappa>0$ such that $K(z)\ge\kappa$ for $|z|<\delta$.
Cover $[0,\ell]$ by finitely many intervals $E\subset[0,\ell]$ of positive length and
diameter less than $\delta$. On each such $E$,
 \[\varphi_\ep(x)\ge\frac{d\kappa}{b(x)+\ep}\int_E\varphi_\ep(y)\dy.\]
Integrating over $E$ and dividing by the positive factor $\int_E\varphi_\ep(y)\dy$, we obtain $\int_E\frac{\dx}{b(x)+\ep}\le\frac1{d\kappa}$. 
Monotone convergence and the finite cover prove \eqref{6.8}.
In particular, $\{b=0\}$ has Lebesgue measure zero.

\textit{Step 2}. For $0<l\le \ell$, define a positive compact operator on $C([0,l])$ by
\begin{equation}\label{6.10}
\mathcal T_l[f](x)=d\int_0^l K(y-x)\frac{f(y)}{b(y)}\dy.
\end{equation}
For every $M>0$ and $f\in C([0,l])$ with $\|f\|_\infty\le M$, \eqref{6.8} implies
 \[\|\mathcal T_l[f]\|_\infty\le dM\|K\|_\infty\int_0^l\frac{\dy}{b(y)}<\infty.\]
Moreover, for $x,x'\in[0,l]$,
  \[|\mathcal T_l[f](x)-\mathcal T_l[f](x')|
    \le dM\sup_{y\in[0,l]}|K(y-x)-K(y-x')|\int_0^l\frac{\dy}{b(y)}. \]
Since $K$ is uniformly continuous on $[-l,l]$, the right-hand side tends to zero as $|x-x'|\to0$, uniformly for $\|f\|_\infty\le M$. Thus the functions $\mathcal T_l[f]$ with $\|f\|_\infty\le M$ are uniformly bounded and equicontinuous. The Arzel\`a--Ascoli theorem implies that $\mathcal T_l$ is compact. 

Since $b>0$ almost everywhere and $K(z)\ge\kappa$ for $|z|<\delta$,  by iterating \eqref{6.10}, we obtain $\mathcal T_l^N[f]\gg0$ for every nonzero $f\in C([0,l])$ with $f\ge0$,  whenever $N\delta>l$. Also $\mathcal T_l[\mathbf1]\gg0$, so $r(\mathcal T_l)>0$. The Krein--Rutman theorem therefore provides a strongly positive eigenfunction and a strictly positive dual eigenfunctional.

\textit{Step 3}. We prove $r(\mathcal T_\ell)\le1$. 
For $l=\ell$, define $\mathcal T_{\ell,\ep}$ by replacing $b$ with
$b+\ep$. Setting $f_\ep=(b+\ep)\varphi_\ep\gg0$ and using \eqref{6.9} we obtain
  \bee\mathcal T_{\ell,\ep}[f_\ep]\le f_\ep,\;\;\Longrightarrow\; r(\mathcal T_{\ell,\ep})\le1.\zzz\label{6.11}\eee
Moreover,
 \[\|\mathcal T_{\ell,\ep}-\mathcal T_\ell\|\le d\|K\|_\infty\int_0^\ell
\left(\frac1b-\frac1{b+\ep}\right)\dx\to0\;\;\;\text{as}\;\;\ep\to 0\]
by dominated convergence. If $r(\mathcal T_\ell)>1$, choose $\varphi\gg0$ such that
$\mathcal T_\ell[\varphi]=r(\mathcal T_\ell)\varphi$ and fix
$1<q<r(\mathcal T_\ell)$. Since $\min_{[0,\ell]}\varphi>0$ and
$\|\mathcal T_{\ell,\ep}-\mathcal T_\ell\|\to0$, for all sufficiently small $\ep>0$, we have
 $\mathcal T_{\ell,\ep}[\varphi]\ge q\varphi$ and $r(\mathcal T_{\ell,\ep})\ge q>1$,
contradicting \eqref{6.11}. Hence $r(\mathcal T_\ell)\le1$.

\textit{Step 4}. For $0<l<\ell$, we claim that
\begin{equation}\label{6.13}
r(\mathcal T_l)<r(\mathcal T_\ell)\le1.
\end{equation}
Let $\varphi_\ell\gg0$ be an eigenfunction associated with $r(\mathcal T_\ell)$. For $x\in[0,l]$, we have
  \[
 r(\mathcal T_\ell)\varphi_\ell(x)-\mathcal T_l[\varphi_\ell|_{[0,l]}](x)=d\int_l^\ell K(y-x)\frac{\varphi_\ell(y)}{b(y)}\dy\ge0.
    \]
The inequality is strict at $x=l$, since $K>0$ near zero, $\varphi_\ell\gg0$, and $b>0$ almost everywhere. Let $\psi_l$ be a strictly positive dual eigenfunctional associated with $r(\mathcal T_l)$. Applying $\psi_l$ to the preceding difference gives
 \[\bigl(r(\mathcal T_\ell)-r(\mathcal T_l)\bigr)\langle\psi_l,\varphi_\ell|_{[0,l]}\rangle>0.\zzz\]
Since $\langle\psi_l,\varphi_\ell|_{[0,l]}\rangle>0$, the preceding inequality implies \eqref{6.13}.

\textit{Step 5}. We now prove \eqref{6.7}. For $\tau>0$, set $V_\tau=\int_0^\tau v(t,\cdot)\dt$ and $k_\tau=bV_\tau$. Since $k_\tau/b=V_\tau$ almost everywhere, $\mathcal T_l[k_\tau]=d\mathcal K_{(0,l)}[V_\tau]$
at every point, including zeros of $b$. Integrating the equation $v_t=\mathcal J_{a,l}[v]$ over $[0,\tau]$ and using $v(0,x)=g(x)$, we obtain, for $x\in[0,l]$,
 \[\begin{aligned}
 v(\tau,x)-g(x) &=d\int_0^\tau\int_0^l K(y-x)v(t,y)\dy\dt-b(x)\int_0^\tau v(t,x)\dt\\
&=d\int_0^l K(y-x)V_\tau(y)\dy-k_\tau(x)\\
 &=\mathcal T_l[k_\tau](x)-k_\tau(x).
    \end{aligned} \]
Rearranging and using $v(\tau,\cdot)\ge0$, we obtain $(\mathrm{Id}-\mathcal T_l)[k_\tau]    =g-v(\tau,\cdot)\le g$. Since $(\mathrm{Id}-\mathcal T_l)^{-1}=\sum_{n\ge0}\mathcal T_l^n$ is positive,
  \[0\le V_\tau(x)\le\frac{[(\mathrm{Id}-\mathcal T_l)^{-1}g](x)}{b(x)}
\;\;\;\text{for\; a.e. }x\in(0,l).\]
The right-hand side is integrable. Integrating over $(0,l)$ and letting $\tau\to\infty$, we obtain \eqref{6.7} with the finite constant
\begin{equation*}
C_l=\int_0^l
\frac{[(\mathrm{Id}-\mathcal T_l)^{-1}\mathbf1](x)}{b(x)}\dx.\zzz
\end{equation*}
The proof is complete.\www
\end{proof}

\begin{remark} Estimate \eqref{6.7} does not require $s(\mathcal J_{a,l})<0$. If $s(\mathcal J_{a,l})=0$, the linear semigroup cannot decay exponentially in operator norm. In the proof of Theorem~\ref{th6.1}, we use the space--time integral in \eqref{6.7} to bound the boundary displacement; see \eqref{6.23a} below.\zzz
\end{remark}

\begin{lem}\label{lem6.2} Fix $T>0$, $\nu>0$, and initial data $(u_0,h_0)$. For $\mu$ in a neighborhood of $\nu$ in $(0,\infty)$, choose $\ell$ larger
than all $h_\mu(T)$ and denote the zero extension of $u_\mu$ to $[0,\ell]$ by $\tilde u_\mu$. There is a constant $C_{T,\nu}$ such that
 \begin{equation}\label{6.16}
\sup_{0\le t\le T}\left(|h_\mu(t)-h_\nu(t)|+
 \|\tilde u_\mu(t,\cdot)-\tilde u_\nu(t,\cdot)\|_{L^\infty(0,\ell)}\right)\le C_{T,\nu}|\mu-\nu|.\zzz
\end{equation}
\end{lem}

\begin{proof}\textit{Step 1: Integral equation for the zero extension}. Let $\bar\mu$ be the supremum of the fixed  neighborhood of $\nu$. Since $h'_\mu(t)\le\mu M_0h_\mu(t)$, we have
 \[h_\mu(t)\le h_0e^{\mu M_0t}\le h_0e^{\bar\mu M_0T},\;\;\; 0\le t\le T.\]
Fix $\ell>h_0e^{\bar\mu M_0T}$. Let $\tilde u_0$ be the zero extension of $u_0$ to $[0,\ell]$, and set $\chi_\mu(t,x)=\mathbf1_{[0,h_\mu(t))}(x)$. By the definition of the zero extension,
  \[\mathcal K_{(0,\ell)}[\tilde u_\mu](s,x)=\int_0^\ell K(y-x)\tilde u_\mu(s,y)\dy=\int_0^{h_\mu(s)}K(y-x)u_\mu(s,y)\dy.
 \]
Fix $t\in(0,T]$. We derive the integral equation separately for the following three cases.
        
\textit{Case 1: $0\le x\le h_0$.} Since $h_\mu$ is strictly increasing, $x<h_\mu(s)$ for every $s\in(0,t]$. Integrating the equation for $u_\mu(s,x)$ from $0$ to $t$ yields
 \[ \tilde u_\mu(t,x)=u_0(x)+
 \int_0^t\left[d\chi_\mu\mathcal K_{(0,\ell)}[\tilde u_\mu](s,x)
+a(x)\tilde u_\mu(s,x)-\tilde u_\mu(s,x)^2\right]\mathrm{d}s.\]
 %Here $u_0(x)=\tilde u_0(x)$ and $\chi_\mu(s,x)=1$ for $s\in(0,t]$.
        
\textit{Case 2: $h_0<x<h_\mu(t)$.} Since $h_\mu$ is continuous and strictly increasing, its inverse $h_\mu^{-1}$ is well defined on $[h_0,h_\mu(T)]$.  In the present case,
 \[ 0<h_\mu^{-1}(x)<t, \;\;\; h_\mu(h_\mu^{-1}(x))=x. \]
The boundary condition therefore implies $u_\mu(h_\mu^{-1}(x),x)=0$. For $s>h_\mu^{-1}(x)$, the point $x$ lies inside the moving interval, so integration of the equation from $h_\mu^{-1}(x)$ to $t$ yields
 \[\tilde u_\mu(t,x)=\int_{h_\mu^{-1}(x)}^t
\left[d\mathcal K_{(0,\ell)}[\tilde u_\mu](s,x)+a(x)\tilde u_\mu(s,x)-\tilde u_\mu(s,x)^2\right]\mathrm{d}s. \]
Since $h_\mu$ is strictly increasing,
  \[ \begin{aligned}
  0\le s\le h_\mu^{-1}(x)
 &\; \Longrightarrow\ h_\mu(s)\le x
 \; \Longrightarrow\ \tilde u_\mu(s,x)=\chi_\mu(s,x)=0,\\
h_\mu^{-1}(x)<s\le t&\; \Longrightarrow\ x<h_\mu(s) \; \Longrightarrow\ \chi_\mu(s,x)=1.
\end{aligned}  \]
        
\textit{Case 3: $h_\mu(t)\le x\le \ell$.} By monotonicity of $h_\mu$, we have $x\ge h_\mu(s)$ for every $s\in[0,t]$. Hence $\tilde u_0(x)=0$ and $\tilde u_\mu(s,x)=\chi_\mu(s,x)=0$ for $0\le s\le t$. 
        
Combining these cases, we obtain
 \begin{equation}\label{6.17}
 \tilde u_\mu(t,x)=\tilde u_0(x)+\int_0^t\left[
d\chi_\mu\mathcal K_{(0,\ell)}[\tilde u_\mu]+a\tilde u_\mu-\tilde u_\mu^2\right](s,x)\mathrm{d}s
        \end{equation}
for every $(t,x)\in[0,T]\times[0,\ell]$, where the identity at $t=0$ follows from the initial condition.
    
\textit{Step 2: An $L^1$ estimate and a boundary estimate.}   Set
 \[ w=\tilde u_\mu-\tilde u_\nu, \;\;\;
 E(t)=\|w(t,\cdot)\|_{L^1(0,\ell)},\;\;\; E_h(t)=|h_\mu(t)-h_\nu(t)|.
 \]
Using the decomposition
 \[\chi_\mu\mathcal K_{(0,\ell)}[\tilde u_\mu]-\chi_\nu\mathcal K_{(0,\ell)}[\tilde u_\nu]=\chi_\mu\mathcal K_{(0,\ell)}[w]
+(\chi_\mu-\chi_\nu)\mathcal K_{(0,\ell)}[\tilde u_\nu],
 \]
together with $\|\mathcal K_{(0,\ell)}[w]\|_1\le\|w\|_1$, $0\le\mathcal K_{(0,\ell)}[\tilde u_\nu]\le M_0$, and $\|\chi_\mu-\chi_\nu\|_1=E_h$, we obtain from \eqref{6.17} that
 \[ E(t)\le \int_0^t\left[ (d+\|a\|_\infty+2M_0)E(s)+dM_0E_h(s) \right]\mathrm{d}s.\www
\]
        
Next, put $Q(r)=\int_r^\infty K(z)\dz$. Then $0\le Q\le1$ and $|Q(r)-Q(s)|\le\|K\|_\infty|r-s|$. Using the zero extensions, we subtract the boundary equations to obtain
 \[\begin{aligned}
 h'_\mu(t)-h'_\nu(t) &=(\mu-\nu)\int_0^\ell \tilde u_\mu(t,x)Q(h_\mu(t)-x)\dx+\nu\int_0^\ell
w(t,x)Q(h_\mu(t)-x)\dx\\
 &\;\;\;+\nu\int_0^\ell \tilde u_\nu(t,x) \bigl[Q(h_\mu(t)-x)-Q(h_\nu(t)-x)\bigr]\dx.
\end{aligned}\]
The bounds for $Q$ and the common initial boundary imply
 \[ E_h(t) \le\int_0^t|h'_\mu(s)-h'_\nu(s)|\,\mathrm{d}s\le M_0\ell t|\mu-\nu|+\nu\int_0^t
\left[E(s)+M_0\ell\|K\|_\infty E_h(s)\right]\mathrm{d}s.
 \]
Adding the two estimates, we deduce from Gronwall's inequality that
 \begin{equation}\label{6.18}
\sup_{0\le t\le T}\bigl(E(t)+E_h(t)\bigr) \le C_1|\mu-\nu|,
        \end{equation}
where $C_1$ is independent of $\mu$ in the fixed neighborhood of $\nu$.\zzz

\textit{Step 3: An $L^\infty$ estimate.} The free boundary equation, positivity of $K$ near zero, and positivity of $u_0$ on $[0,h_0)$ imply $h'_\nu(0)>0$. Since $h'_\nu$ is continuous and positive on $[0,T]$, set
    \[c=\min_{0\le t\le T}h'_\nu(t)>0, \;\;\;\Delta=\sup_{0\le t\le T}E_h(t). \]
For fixed $x\in[0,\ell]$, if $\chi_\mu(s,x)\ne\chi_\nu(s,x)$, then $|h_\nu(s)-x|\le\Delta$. For any two such times $s_1<s_2$, we have $c(s_2-s_1)\le h_\nu(s_2)-h_\nu(s_1)\le 2\Delta$. Thus the set of times at which the indicators differ has measure at most $2\Delta/c$, and hence, $\int_0^T|\chi_\mu(s,x)-\chi_\nu(s,x)|\,\mathrm{d}s\le\frac{2\Delta}{c}$.
    
Set $F(t)=\|\tilde u_\mu(t,\cdot)-\tilde u_\nu(t,\cdot)\|_\infty$ and $C=\|a\|_\infty+2M_0$.
Subtracting \eqref{6.17} for $\mu$ and $\nu$ and using $\|\mathcal K_{(0,\ell)}[w]\|_\infty\le\|K\|_\infty E$, we obtain
 \[|\tilde u_\mu(t,x)-\tilde u_\nu(t,x)|\le d\|K\|_\infty\int_0^t E(s)\,\mathrm{d}s
 +dM_0\int_0^t |\chi_\mu(s,x)-\chi_\nu(s,x)|\,\mathrm{d}s+C\int_0^t F(s)\,\mathrm{d}s.\]
Taking the supremum over $x\in[0,\ell]$ yields 
 \begin{equation}\label{6.19}
  F(t)\le d\|K\|_\infty\int_0^t E(s)\,\mathrm{d}s+\frac{2dM_0\Delta}{c} +C\int_0^t F(s)\,\mathrm{d}s.
    \end{equation}
By \eqref{6.18}, both $E(s)$ and $\Delta$ are bounded by $C_1|\mu-\nu|$. The conclusion \eqref{6.16} now follows from \eqref{6.19} and Gronwall's inequality.\www
\end{proof}

\begin{proof}[Proof of Theorem~\ref{th6.1}] \textup{(i)} Fix $\mu>0$ and write $\bar h=h_\mu(\infty)<\infty$. Suppose that $s(\mathcal J_{a,\bar h})>0$.
Choose $\sigma,\kappa>0$ such that $K(z)\ge\kappa$ for $|z|<\sigma$. By spectral continuity, choose
$0<\ep<\frac13\min\{\sigma,\bar h\}$ with $s(\mathcal J_{a,\bar h-\ep})>0$.
For all sufficiently large $t$, $h_\mu(t)>\bar h-\ep$. Comparison with the logistic equation on $[0,\bar h-\ep]$ and convergence to its positive equilibrium imply
\begin{equation}\label{6.20}
\liminf_{t\to\infty}\min_{x\in[0,\,\bar h-\ep]}
\bigl(u_\mu(t,x)-U_{\bar h-\ep}(x)\bigr)\ge0.
\end{equation}
Restricting the boundary flux to $x\in[\bar h-2\ep,\bar h-\ep]$ and $y\in[\bar h,\bar h+\ep]$, we obtain
  \[ h'_\mu(t)\ge\mu\int_{\bar h-2\ep}^{\bar h-\ep}\int_{\bar h}^{\bar h+\ep}
K(y-x)u_\mu(t,x)\dy\dx\ge\mu\kappa\ep\int_{\bar h-2\ep}^{\bar h-\ep}u_\mu(t,x)\dx.\]
This together with \eqref{6.20} yields
  \[\liminf_{t\to\infty}h'_\mu(t)\ge\mu\kappa\ep\int_{\bar h-2\ep}^{\bar h-\ep}U_{\bar h-\ep}(x)\dx>0,\]
contrary to $h_\mu(t)\to \bar h$.
Thus $s(\mathcal J_{a,\bar h})\le0$ and $\bar h\le l_*$ by the definition \eqref{6.5a} for $l_*$.
The fixed-interval logistic solution on $[0,\bar h]$ starting from the
zero extension of $u_0$ dominates $u_\mu$ and converges uniformly to
zero by \cite[Theorem B]{WZ26JMPA}. This proves part~\textup{(i)}.

\textup{(ii)} Since $h_\mu(t)>h_0$ for $t>0$, a finite limit is impossible when
$h_0\ge l_*$. Part~\textup{(ii)} follows.

\textup{(iii)} Assume $h_0<l_*$. We first show that vanishing happens for small $\mu$. Fix $h_0<l<l_*$.
Let $g$ be the zero extension of $u_0$ to $[0,l]$, and let
$v_t=\mathcal J_{a,l}[v]$, $v(0)=g$.
Let $C_l$ be given by Lemma \ref{lem6.1} and choose $\ud \mu>0$ such that
 \bee
 h_0+\ud\mu C_l\|u_0\|_\infty<l.\label{6.22a}\eee
We claim that $h_\mu(t)<l$ for all $t>0$ and all $0<\mu\le\ud\mu$, which implies $h_\mu(\infty)<\infty$, and $\lim_{t\to\infty}\|u_\mu(t,\cdot)\|_\infty=0$ for all $0<\mu\le\ud\mu$. Assume on the contrary that  $h_\mu(T_{\mu,l})=l$ for some $T_{\mu,l}>0$. Comparison gives $u_\mu\le v$ on $[0,T_{\mu,l}]\times[0,l]$. Since $s(\mathcal J_{a,l_*})=0$, we may take $\ell=l_*$ in Lemma~\ref{lem6.1}. Integrating the equation for the free boundary yields 
 \bee 
 l-h_0\le\mu\int_0^{T_{\mu,l}}\!\!\int_0^{h_\mu(s)}\!\!u_\mu(s,x)\dx\mathrm{d}s
\le\mu\int_0^{T_{\mu,l}}\!\!\int_0^l\!\! v(s,x)\dx\mathrm{d}s
\le\mu C_l\|u_0\|_\infty.\label{6.23a}\eee
This contradicts \eqref{6.22a}. 

We next show that spreading happens for large $\mu$. For every fixed $T>0$, we claim that
\begin{equation}\label{6.23}
\lim_{\mu\to\infty}h_\mu(T)=\infty.
\end{equation}
Otherwise, monotonicity in $\mu$ implies $h_\mu(T)\le \ell$ for some $\ell$ and all $\mu$. Extend $u_\mu(t,\cdot)$ by zero to $[0,\ell]$ and denote the monotone pointwise limit of $u_\mu(t,x)$ by $\bar u(t,x)$; write
$\bar h(t)=\lim_{\mu\to\infty}h_\mu(t)$. For $t>0$ and $x<\bar h(t)$, some finite $\mu$ satisfies $x<h_\mu(t)$, so $\bar u(t,x)>0$. The integrated equation for the free boundary, $u_\mu\le M_0$, and dominated convergence then yield
\[
0=\lim_{\mu\to\infty}\frac{h_\mu(T)-h_0}{\mu}
=\int_0^T\int_0^{\bar h(t)}\int_{\bar h(t)}^\infty
K(y-x)\bar u(t,x)\dy\dx\dt>0.
\]
The last integral is positive because $K(0)>0$ and $\bar u(t,x)>0$ for $0\le x<\bar h(t)$. This contradiction proves \eqref{6.23}. Thus $h_\mu(T)>l_*$ for large $\mu$,
and part~\textup{(i)} implies spreading.

Now we determine the threshold parameter and the corresponding limiting boundary. Set
  \[\Sigma=\{\mu>0:h_\mu(\infty)<\infty\},\;\;\;\mu^*=\sup\Sigma.\]
Then $h_\mu(\infty)=\infty$ for all $\mu>\mu^*$. Comparison with respect to $\mu$ shows that if $\mu\in\Sigma$, then $(0,\mu]\subset\Sigma$. The above arguments imply $0<\mu^*<\infty$.
If $h_{\mu^*}(\infty)=\infty$, then $h_{\mu^*}(T)>l_*$ for some finite $T$. Lemma~\ref{lem6.2} ensures that $h_\mu(T)>l_*$ for a parameter $\mu<\mu^*$ sufficiently close to $\mu^*$. This contradicts part~\textup{(i)} due to $(0,\mu^*)\subset\Sigma$. Therefore
$\Sigma=(0,\mu^*]$, and vanishing occurs at $\mu^*$.

Finally, part~\textup{(i)} and vanishing at $\mu^*$ imply $h_{\mu^*}(\infty)\le l_*$. Suppose that $h_{\mu^*}(\infty)<l_*$ and choose $h_{\mu^*}(\infty)<l<l_*$. Let $C_l$ be given in Lemma~\ref{lem6.1} with $\ell=l_*$. Since $\|u_{\mu^*}(t,\cdot)\|_{L^\infty(0,l)}\to 0$ as $t\to\infty$, we can choose $T>0$ such that
 \[  h_{\mu^*}(T)+\mu^*C_l\|u_{\mu^*}(T,\cdot)\|_{L^\infty(0,l)}<l. \]
By the continuity (Lemma~\ref{lem6.2}), there exists $\mu>\mu^*$ sufficiently close to $\mu^*$ such that
  \begin{equation}\label{6.24}
h_\mu(T)+\mu C_l\|u_\mu(T,\cdot)\|_{L^\infty(0,l)}<l.
    \end{equation}
For such a fixed $\mu$, we have $h_\mu(\infty)=\infty$. There exists $T_\mu>T$ such that $h_\mu(T_\mu)=l$. Similar to the derivation of \eqref{6.23a}, 
  \[ l-h_\mu(T)\le\mu\int_T^{T_\mu}\!\!\int_0^{h_\mu(s)}\!
u_\mu(s,x)\dx\,\mathrm{d}s\le\mu C_l\|u_\mu(T,\cdot)\|_{L^\infty(0,l)}. \]
This contradicts \eqref{6.24}. Thus $h_{\mu^*}(\infty)=l_*$, proving \eqref{6.6}.\www
\end{proof}

The following lemma presents a sufficient condition to guarantee \eqref{6.5}.\zzz

\begin{lemma} \label{lem6.3} For $l>0$, define
\begin{equation*}
\zeta_K(l)=\frac12\int_{\mathbb R}K(z)
\min\left\{\frac{\pi^2z^2}{l^2},4\right\}\dz.
\end{equation*}
Then $\zeta_K(l)\to0$ as $l\to\infty$. If $m\ge m_0$ on an interval $I=(r,r+\ell)\subset(0,\infty)$, then
\begin{equation}\label{6.25}
s(\mathcal J_{a,I})\ge m_0-d\zeta_K(\ell).\zzz
\end{equation}
In particular, $m_0>d\zeta_K(\ell)$ implies \eqref{6.5}.\zzz
\end{lemma}

\begin{proof} First, $0\le K(z)\min\big\{\frac{\pi^2z^2}{l^2},4\big\}\le4K(z)$. Since $K\in L^1(\mathbb R)$ and the integrand converges pointwise to zero as $l\to\infty$, the dominated convergence theorem gives $\zeta_K(l)\to0$.

Extend $\varphi(x)=\sin(\pi(x-r)/\ell)$ from $I$ by zero to $\mathbb R$, and write $\|\cdot\|_2=\|\cdot\|_{L^2(\mathbb R)}$. The Rayleigh formula in Theorem~\ref{th2.1}\textup{(v)}  yields
\[
s(\mathcal J_{a,I})\ge\frac{\int_I m\varphi^2\dx
-\frac d2\int_{\mathbb R}K(z)\|\varphi(\cdot+z)-\varphi\|_2^2\dz}
{\|\varphi\|_2^2}.\]
Since the zero extension belongs to $H^1(\mathbb R)$,
the fundamental theorem of calculus and the Cauchy--Schwarz
inequality yield
\[
\begin{aligned}
\|\varphi(\cdot+z)-\varphi\|_2^2=z^2\int_{\mathbb R}
\left|\int_0^1\varphi'(x+\theta z)\,\mathrm{d}\theta\right|^2\dx
\le z^2\|\varphi'\|_2^2.\end{aligned}
\]
The triangle inequality and direct integration also give
 \[\|\varphi(\cdot+z)-\varphi\|_2^2\le4\|\varphi\|_2^2,
 \;\;\; \|\varphi'\|_2^2/\|\varphi\|_2^2=\pi^2/\ell^2.\]
Combining these estimates with $m\ge m_0$ on $I$, we obtain
\[
s(\mathcal J_{a,I})\ge m_0-\frac d2\int_{\mathbb R}K(z)
\min\left\{\frac{\pi^2z^2}{\ell^2},4\right\}\dz=m_0-d\zeta_K(\ell),\]
which proves \eqref{6.25}.

If $m_0>d\zeta_K(\ell)$, choose $l>r+\ell$. Domain monotonicity implies
$s(\mathcal J_{a,l})\ge s(\mathcal J_{a,I})>0$, proving \eqref{6.5}.\www
\end{proof}

The following conditions are nonlocal counterparts of those used
in \cite{W15JDE}.\zzz

\begin{prop}\label{prop6.1}
Condition \eqref{6.5} holds for every $d>0$ if one of the
following conditions is satisfied.\www
\begin{enumerate}[label=\textup{(\roman*)}]
\item There exist a constant $\rho>0$ and intervals
$[x_i,y_i]\subset[0,\infty)$ such that
 \[y_i-x_i\to\infty,\;\;\; m(x)\ge\rho
\;\;\;\text{for }x\in[x_i,y_i],\;\;\; i\ge1.\zzz\]
\item For some $p\in(0,2]$,
\[M_p:=\int_{\mathbb R}|z|^pK(z)\dz<\infty.\]
There exist $\rho>0$, $k>1$, $x_i\to\infty$ and $-p<\gamma\le0$
such that $m(x)\ge\rho x^\gamma$ on $[x_i,kx_i]$.\www
\end{enumerate}
\end{prop}

\begin{proof} Part (i) follows from $\zeta_K(y_i-x_i)\to0$ as $i\to\infty$. For part (ii), the inequality $\min\{r,4\}\le r^{p/2}4^{1-p/2}$ for $r\ge0$ implies
 \[\zeta_K\big((k-1)x_i\big)\le 2^{1-p}\pi^p M_p[(k-1)x_i]^{-p}.\]
Since $\gamma\le0$, we have $m(x)\ge\rho k^\gamma x_i^\gamma$ on
$[x_i,kx_i]$. Lemma~\ref{lem6.3} therefore yields
 \[s\big(\mathcal J_{a,(x_i,kx_i)}\big)\ge\rho k^\gamma x_i^\gamma-d\,2^{1-p}\pi^pM_p(k-1)^{-p}x_i^{-p}>0\]
for large $i$, because $p+\gamma>0$. Domain monotonicity then implies \eqref{6.5}.
\end{proof}

\subsection{Positive equilibria and long-time behavior}

The next comparison also holds when the larger interval is unbounded.

\begin{lem}\label{lem6.4}
Let $I_1\subset I_2\subset(0,\infty)$ be intervals, with $I_1$ bounded. Assume $s(\mathcal J_{a,I_1})>0$ and let $U_1\in C(\bar I_1)$ be the unique positive solution of
  \[d\int_{I_1}K(y-x)U_1(y)\dy+(m(x)-d)U_1(x)-U_1(x)^2=0,\;\;\; x\in\bar I_1.\]
If $U_2\in C(\bar I_2)$ is positive on $\bar I_2$ and satisfies
  \[d\int_{I_2}K(y-x)U_2(y)\dy+(m(x)-d)U_2(x)-U_2(x)^2=0,\;\;\; x\in\bar I_2,\]
then $U_1\le U_2$ on $\bar I_1$.
\end{lem}

\begin{proof} Restricting $U_2$ to $\bar I_1$ gives
\begin{equation}\label{6.28}
\mathcal J_{a,I_1}[U_2](x)-U_2(x)^2
=-d\int_{I_2\setminus I_1}K(y-x)U_2(y)\dy\le0,
\;\;\; x\in\bar I_1.
\end{equation}
The continuous positive function $U_2$ has a positive minimum on $\bar I_1$.
Set $\theta=\max_{\bar I_1}U_1/U_2$.
If $\theta>1$, then $W=\theta U_2-U_1\ge0$ has a zero at some
$x_*\in\bar I_1$. By \eqref{6.28} and the equation for $U_1$,
\[
\begin{aligned}
0\le d\int_{I_1}K(y-x_*)W(y)\dy
=\mathcal J_{a,I_1}[W](x_*)
\le\theta U_2(x_*)^2-U_1(x_*)^2
=-\theta(\theta-1)U_2(x_*)^2<0,
\end{aligned}
\]
a contradiction. Thus $\theta\le1$, which proves $U_1\le U_2$.
\end{proof}

\begin{prop}\label{prop6.3}
Assume \eqref{6.1}, \eqref{6.2}, and \eqref{6.5}. For $l>l_*$, let $U_l$ denote the unique positive solution to \eqref{6.4}. Then $U_l\to U_*$ locally uniformly as $l\to\infty$, where $U_*$ is the minimal bounded continuous positive solution to
\begin{equation}\label{6.29}
d\int_0^\infty K(y-x)U(y)\dy+(m(x)-d)U-U^2=0,\;\;\; x\ge0.
\end{equation}
Moreover, if $h_\mu(\infty)=\infty$, then for every $R>0$,
 \bee
\liminf_{t\to\infty}\min_{x\in[0,R]}\bigl(u_\mu(t,x)-U_*(x)\bigr)\ge0.
 \label{6.30}\eee
In particular, spreading implies persistence on every fixed compact interval. \zzz
\end{prop}

\begin{proof} By Lemma~\ref{lem6.4}, the equilibria $U_l$, $l>l_*$, increase with $l$. The constant $M=\max\{\sup m,0\}$ is an upper solution, so $0<U_l\le M$. Thus $U_l(x)\nearrow U_*(x)$ for each $x\ge0$. For any fixed $x$, choose $l>\max\{l_*,x\}$; then $U_*(x)\ge U_l(x)>0$. 
Dominated convergence in the equilibrium equations shows that $U_*$ 
satisfies \eqref{6.29}. Solving \eqref{6.29} for its positive root, we obtain
  \[U_*(x)=\frac{a(x)+\sqrt{a(x)^2+4d\int_0^\infty K(y-x)U_*(y)\dy}}2.\]
 Since $\int_0^\infty K(y-x)U_*(y)\dy$ is continuous, we conclude that $U_*$ is continuous and Dini's theorem ensures local uniform convergence of $U_l$ to $U_*$.
If $U$ is any bounded continuous positive solution of \eqref{6.29}, Lemma~\ref{lem6.4} implies $U_l\le U$ on $[0,l]$. Letting $l\to\infty$ proves $U_*\le U$.

Assume $h_\mu(\infty)=\infty$ and fix $R>0$.
For each $l>\max\{l_*,R\}$, choose $T_l$ with $h_\mu(T_l)>l$.
Comparison from $T_l$ onward with the logistic equation on $[0,l]$,
whose positive solutions converge to $U_l$, yields $\liminf_{t\to\infty}\min_{x\in[0,R]}
\bigl(u_\mu(t,x)-U_l(x)\bigr)\ge0$. 
For any given $\ep>0$, there exists $l$ such that $\|U_l-U_*\|_{C([0,R])}<\ep$. It then follows that 
  \[\liminf_{t\to\infty}\min_{x\in[0,R]}\bigl(u_\mu(t,x)-U_*(x)\bigr)\ge-\ep.\]
Letting $\ep\searrow0$ proves \eqref{6.30}.\www
\end{proof}

\begin{theo}
\label{th6.2}
Assume \eqref{6.1}--\eqref{6.2} and
\begin{equation}\label{6.32}
0<m_1:=\liminf_{x\to\infty}m(x)\le \limsup_{x\to\infty}m(x)=:m_2<\infty.
 \end{equation}
Then \eqref{6.5} holds, \eqref{6.29} has a unique bounded continuous positive solution $U$, and $\inf_{x\ge0}U(x)>0$. Moreover,
 \bee
 m_1\le\liminf_{x\to\infty}U(x)\le\limsup_{x\to\infty}U(x)\le m_2.\zzz\label{6.33}\eee
If $h_\mu(\infty)=\infty$, then 
   \bee
 u_\mu(t,\cdot)\to U(\cdot)\;\;\;\text{locally uniformly on }[0,\infty)
 \;\;\;(t\to\infty).
 \label{6.34}\eee
\end{theo}

\begin{proof} Proposition~\ref{prop6.1}(i) implies \eqref{6.5}, and
Proposition~\ref{prop6.3} provides the minimal bounded continuous
positive equilibrium $U_*$ of \eqref{6.29}.

\textit{Step 1: The estimates \eqref{6.33}}. Let $U$ be any bounded continuous positive solution of \eqref{6.29}. Then $U\le\sup m$. By \eqref{6.32}, for any given $0<\varepsilon<m_1$, there exists $l\gg1$ such that $m_1-\varepsilon<m(x)<m_2+\varepsilon$ for all $x\ge l$. It is clear that the problem
\[
d\int_l^\infty K(y-x)U(y)\dy+(m_1-\varepsilon-d)U-U^2=0,\;\;\; x\ge l
\]
has a unique bounded positive solution $\widetilde U\in C([l,\infty))$, and $\lim_{x\to\infty}\widetilde U(x)=m_1-\varepsilon$. The comparison principle gives $U\ge \widetilde U$ on $[l,\infty)$. The first estimate in \eqref{6.33} follows, as $\varepsilon>0$ is arbitrary.

Let $x_i\to\infty$ satisfy $\lim_{i\to\infty}U(x_i)=\limsup_{x\to\infty}U(x)=:\gamma$. 
Then $m_1\le\gamma\le\sup m$, and $\limsup_{i\to\infty}U(x_i+z)\le\gamma$ for each fixed $z\in\mathbb{R}$. Moreover,
 \begin{align}
0&=d\int_0^\infty K(y-x_i)U(y)\dy+\bigl(m(x_i)-d\bigr)U(x_i)-U^2(x_i)\nonumber\\
&\le d\int_{-x_i}^\infty K(z)U(x_i+z)\dz+\bigl(m_2+\varepsilon-d\bigr)U(x_i)-U^2(x_i)\nonumber\\
&=d\int_\mathbb{R} K(z)U(x_i+z)\mathbf{1}_{\{z\ge -x_i\}}\dz
+\bigl(m_2+\varepsilon-d\bigr)U(x_i)-U^2(x_i).\label{6.35}
 \end{align}
As $\int_{\mathbb R}K(z)\dz=1$ and $0<U\le\sup m$, for any $\delta>0$, there exists a large $l>0$ such that 
 \[d\int_{\mathbb{R}\setminus(-l,l)} K(z)U(x_i+z)\mathbf{1}_{\{z\ge -x_i\}}\dz<\delta,\;\;\;\forall\;i.\]
For such a fixed $l$, there exists a sufficiently large $i_0$ such that $U(x_i+z)\le\gamma+\delta$ for all $|z|\le l$ and all $i\ge i_0$. It follows that, for all $i\ge i_0$, 
 \begin{align*}
 d\int_\mathbb{R} K(z)U(x_i+z)\mathbf{1}_{\{z\ge -x_i\}}\dz&=
 d\left(\int_{\mathbb{R}\setminus(-l,l)}+\int_{-l}^l\right) K(z)U(x_i+z)\mathbf{1}_{\{z\ge -x_i\}}\dz\\
 &\le\delta+d(\gamma+\delta).
 \end{align*}
Letting $i\to\infty$ in \eqref{6.35} we obtain $0\le \delta+d \delta+(m_2+\varepsilon)\gamma-\gamma^2$. The arbitrariness of $\varepsilon>0$ and $\delta>0$  implies the last estimate in \eqref{6.33}.

\textit{Step 2: Uniqueness.} Minimality implies $U_*\le U$. By Step~1, $\liminf_{x\to\infty}U_*(x)\ge m_1>0$. Together with the continuity and positivity of $U_*$ on compact intervals, this implies $\inf_{x\ge0}U_*(x)>0$. By Step~1 and boundedness,  $1\le k=\sup_{x\ge0}U(x)/U_*(x)<\infty$. If $k>1$, choose $x_i\ge0$ with $U(x_i)/U_*(x_i)\to k$ and put $W=kU_*-U\ge0$. Then $W(x_i)\to0$, and subtraction of the equilibrium equations gives\zzz
 \[\begin{aligned}
d\int_0^\infty K(y-x_i)W(y)\dy&=kU_*^2(x_i)-U^2(x_i)-(m(x_i)-d)W(x_i)\\
&=-k(k-1)U_*^2(x_i)+(kU_*(x_i)+U(x_i)-m(x_i)+d)W(x_i).
\end{aligned}\]
This is a contradiction because $\inf U_*>0$ and $W(x_i)\to0$. Hence $k=1$, and therefore $U=U_*$.

\textit{Step 3: The proof of \eqref{6.34}}. Suppose $h_\mu(\infty)=\infty$, and let $v$ be the unique solution of \zzz
 \[v_t=d\int_0^\infty K(y-x)v(t,y)\dy+(m(x)-d)v-v^2,\;\;\; v(0,x)=M_1,\zzz\]
where $M_1=\|u_0\|_\infty+2\sup m$. By comparison, $v\ge u_\mu$, and by \eqref{6.30}, \zzz
 \bee
 \liminf_{t\to\infty}v(t,\cdot)\ge\liminf_{t\to\infty}u_\mu(t,\cdot)\ge U(\cdot)\;\;\;\text{locally uniformly on}\;\;[0,\infty). 
 \label{6.36}\eee
As in the proof of \cite[Lemma 2.8]{LLW22JDE}, we can prove that $v(t,\cdot)\to U(\cdot)$ locally uniformly on $[0,\infty)$ as $t\to\infty$. Combining this with \eqref{6.36}, we obtain \eqref{6.34}.\www
\end{proof}

Theorem~\ref{th6.2} does not require $m(x)$ to converge at infinity.
The condition $\liminf_{x\to\infty}m(x)>0$ ensures
$\inf_{x\ge0}U_*(x)>0$, which is used in the uniqueness proof.
When $m(x)>0$ and $m(x)\to0$ as $x\to\infty$,
our arguments do not establish uniqueness of a bounded positive
equilibrium or convergence of spreading solutions.
\vskip 4pt

 {\noindent \bf Ethical Statement.}

\begin{itemize}
    \item Originality: The work presented in this manuscript is original and has not been published elsewhere in any form or language. 
    \item Simultaneous Submission: This manuscript has not been submitted to more than one journal for simultaneous consideration.
    \item Availability of data and material: Not applicable.
    \item Conflict of Interest: There are no conflicts of interest related to this manuscript. 
\end{itemize}

 \end{document}